\documentclass[A4paper,12pt]{book}
\usepackage{amssymb,amsmath}
\usepackage[russian]{babel}
\usepackage{graphicx}
\newcount\exnom\exnom=0

\long\def\exo/{\vspace{0.01cm} \noindent\advance\exnom by1{\bf
{\the\exnom}.}}

\newcommand{\ds}{\displaystyle}

\def\im{\mathop{\rm Im}\nolimits}

\begin{document}
\thispagestyle{empty}

\begin{center} \large \bf JORDAN DERIVATIONS IN FINITE ENDOMORPHISM SEMIRINGS
\end{center}

\begin{center}  Dimitrinka I. Vladeva
\end{center}

\begin{quote}\small
\emph{{Abstract}}
The aim of this article is to start a study of Jordan derivations in finite endomorphism semirings.



\end{quote}

{\large\textbf{1. Introduction and preliminaries.}}
\vspace{2mm}

 The differential algebra has been studied by many authors for the last seventy   years and especially the relationships between
derivations and the structure of rings. The notion of the ring with derivation is  old and plays an important role in the integration of analysis, algebraic geometry and algebra.
In 1950 J. Ritt, [9], and in 1973 E. Kolchin, [7], wrote the  classical  books on differential algebra.

 During the last few decades there has been a great deal of works concerning derivations in rings, in Lie rings, in skew polynomial rings and other algebraic structures.

  Herstein, see [5] and [6], constructed, starting from the ring $R$, a new ring, namely the Jordan ring R, defining the product in this $a\circ b = ab + ba$ for any $a, b \in R$. This new product is well-defined and it can be easily verified that $(R,+,\circ)$ is a ring. An additive mapping $D$, from
the Jordan ring into itself, is said by Herstein to be a Jordan derivation, if $D(a\circ b) = D(a)\circ b + a\circ D(b)$, for
any $a, b \in R$. So, in 1957, Herstein proved a classical result:
\vspace{1mm}

\textsl{If $R$ is a prime ring of a characteristic different from 2, then every
Jordan derivation of $R$ is a derivation.}
\vspace{1mm}

Later in 1988 M. Bre\v{s}ar, see [1] and [2], extended this theorem.
\vspace{3mm}

For  semilattice $\mathcal{M}$ the set $\mathcal{E}_\mathcal{M}$ of the endomorphisms of $\mathcal{M}$ is a semiring
 with respect to the addition and multiplication defined with:

 $\bullet \; h = f + g \; \mbox{when} \; h(x) = f(x)\vee g(x) \; \mbox{for all} \; x \in \mathcal{M}$,

 $\bullet \; h = f\cdot g \; \mbox{when} \; h(x) = f\left(g(x)\right) \; \mbox{for all} \; x \in \mathcal{M}$.

 This semiring is called the \emph{ endomorphism semirimg} of the semilattice $\mathcal{M}$.
In this paper all semilattices considered are finite chains.
 We fix a finite chain\break $\mathcal{C}_n = \; \left(\{0, 1, \ldots, n - 1\}\,,\,\vee\right)\;$ and denote the endomorphism semiring of this chain by ${\mathcal{E}}_{\mathcal{C}_n}$. We do not assume that $\alpha(0) = 0$ for arbitrary $\alpha \in {\mathcal{E}}_{\mathcal{C}_n}$. So, there is not a zero in  endomorphism semiring ${\mathcal{E}}_{\mathcal{C}_n}$.

 Let us fix  elements $a_0,  \ldots, a_{k-1} \in \mathcal{C}_n$, where $k \leq n$ and $a_0 <  \ldots < a_{k-1}$. Let us cosider $A = \{a_0, \ldots, a_{k-1}\}$.
We  shall be interested  endomorphisms $\alpha \in {\mathcal{E}}_{\mathcal{C}_n}$ such that $\im(\alpha) \subseteq A$ and denote this set by
 $\sigma^{(n)}\{a_0,  \ldots, a_{k-1}\}$.

\vspace{1mm}

Let $\{b_0, \ldots, b_{\ell-1}\} \subseteq \{a_0, \ldots, a_{k-1}\}$ and consider the set
 $$\sigma^{(n)}\{b_0, \ldots, b_{\ell-1}\} = \{\beta |\; \beta \in \sigma^{(n)}\{a_0, \ldots, a_{k-1}\}, \im(\beta) = \{b_0, \ldots, b_{\ell-1}\}\; \}$$

For  $\beta_1, \beta_2 \in \sigma^{(n)}\{b_0, \ldots, b_{\ell-1}\}$ let $\beta_1 \sim \beta_2$ if and only if the sets $\im(\beta_1)$ and $\im(\beta_2)$ have a common least element. In this way we define an equivalence relation. Any equivalence class can be identified with its least element which is the constant endomorphism $\overline{b_m}$, where $m = 0, \ldots \ell - 1$.
\vspace{1mm}

Now take a simplicial complex $\Delta$ with vertex set $V = \{\overline{a_0}, \ldots, \overline{a_{k-1}}\}$. The set $\{\overline{b_0}, \ldots, \overline{b_{\ell-1}}\}$ is a subset of $\Delta$. Hence, we can consider the set $\sigma^{(n)}\{b_0, \ldots, b_{\ell-1}\}$ as a face of $\Delta$. In particular, when the simplicial complex $\Delta$ consists of all subsets of $V$, it is called a simplex (see [13]) and  $\Delta = \sigma^{(n)}\{a_0, \ldots, a_{k-1}\}$.

\vspace{2mm}

It is easy to see that
for any set $A = \{a_0,  \ldots, a_{k-1}\} \subseteq \mathcal{C}_n$  the simplex
 $\sigma^{(n)}\{a_0, \ldots, a_{k-1}\}$ is a subsemiring of  ${\mathcal{E}}_{\mathcal{C}_n}$.
The number $k$ is called a {\emph{dimension}} of simplex $\sigma^{(n)}\{a_0, \ldots, a_{k-1}\}$.
Any  simplex $\sigma^{(n)}\{b_0, b_1, \ldots, b_{\ell - 1}\}$, where $b_0, \ldots, b_{\ell-1} \in A$, is  a {\emph{face}} of  simplex $\sigma^{(n)}\{a_0, \ldots, a_{k-1}\}$. If $\ell < k$,  face $\sigma^{(n)}\{b_0, b_1, \ldots, b_{\ell - 1}\}$ is called a {\emph{proper face}}.

\vspace{2mm}

 The proper faces of simplex $\sigma^{(n)}\{a_0, \ldots, a_{k-1}\}$ are:
\vspace{1mm}

$\bullet$ $0$ -- simplices, which are \emph{vertices} $\overline{a_0}, \ldots, \overline{a_k}$.

\vspace{1mm}

$\bullet$ $1$ -- simplices, which are called {\emph{strings}}. They are denoted by $\mathcal{STR}^{(n)}\{a,b\}$, where $a, b \in A$.

\vspace{1mm}

$\bullet$ $2$ -- simplices, which are called {\emph{triangles}}. They are denoted by ${\triangle}^{(n)}\{a,b,c\}$, where $a, b, c \in A$.

\vspace{1mm}

The  endomorphisms $\alpha \in \sigma^{(n)}\{a_0, a_1, \ldots, a_{k-1}\}$ such that
$$\alpha(0) = \cdots = \alpha(i_0-1) = a_0, \alpha(i_0) = \cdots = \alpha(i_0+ i_1 -1) = a_1, \; \cdots$$ $$\alpha(i_0 + \cdots +i_{k-2}) = \cdots = \alpha(i_0 + \cdots + i_{k-1} - 1) = a_{k-1}$$
 we denote  by $\alpha = (a_0)_{i_0}(a_1)_{i_1} \ldots (a_{k-1})_{i_{k-1}}$, where
 $\ds \sum_{p=0}^{k-1} i_p = n$.
\vspace{1mm}

 The endomorphism semirings of a finite semilattice are well-established, see  [13], [14] and [15]. Basic facts for semirings  can be found in [4]. The results for derivations in semirings had been proved earlier -- [11], [16], [17], [18], [19] and [20]. Concerning background of  simplicial complexes and combinatorics a reader is  referred to [8] and [10].
 \vspace{7mm}

{\large \textbf{2. Ten types Jordan derivations in a triangle}}
\vspace{4mm}

Our  aim of this section is to prove that for any endomorphism $\alpha \in \triangle^{(n)}\{a,b,c\}$ the map
$$\partial_\alpha : \triangle^{(n)}\{a,b,c\} \rightarrow \triangle^{(n)}\{a,b,c\}$$
which is Jordan multiplication, i. e. $\partial_\alpha(\beta) = \alpha\beta + \beta\alpha$, where $\beta \in \triangle^{(n)}\{a,b,c\}$ is a derivation in the triangle and to find to find the maximal subsemiring of $\triangle^{(n)}\{a,b,c\}$ closed under this derivation.
\vspace{1mm}

For any endomorphisms $\alpha, \beta, \gamma \in \triangle^{(n)}\{a,b,c\}$ it follows
$$\partial_\alpha(\beta + \gamma) = \alpha(\beta + \gamma) + (\beta + \gamma)\alpha = \alpha\beta + \beta\alpha + \alpha\gamma + \gamma\alpha = \partial_\alpha(\beta) + \partial_\alpha(\gamma).$$

Hence the map $\partial_\alpha$ for all $\alpha \in \triangle^{(n)}\{a,b,c\}$ is a linear.
\vspace{1mm}

From [15] we know that any endomorphism $\alpha \in \triangle^{(n)}\{a,b,c\}$  can be characterized by ordered triple $(x,y,z)$, where $\alpha(a) = x$, $\alpha(b) = y$ and $\alpha(c) = z$ and $x, y, z \in \{a,b,c\}$ and this triple is called a {\emph{type}} of $\alpha$. This is denoted by $\alpha \in (x,y,z)$. So, there are ten different types in a triangle.

The proof that $\partial_\alpha$ is a derivation for some $\alpha \in \triangle^{(n)}\{a,b,c\}$ depends of the type of this endomorphism. Thus we show ten different proofs such that for any endomorphism $\alpha$ from a given type the Jordan multiplication $\partial_\alpha(\beta) = \alpha\beta + \beta\alpha$, where $\beta \in \triangle^{(n)}\{a,b,c\}$ is a derivation.
\vspace{1mm}

An arbitrary Jordan  multiplication $\partial_\alpha$, where $\alpha \in (x,y,z)$
is denoted by $\partial_{(x,y,z)}$.

\vspace{3mm}

In the first three cases we suppose that $\beta$ and $\gamma$ are arbitrary endomorphisms of $\triangle^{(n)}\{a,b,c\}$.
\vspace{5mm}

\emph{Case 1.} Let $\alpha \in (a,b,c)$. Following [15] this means that $\alpha$ is right identity of semiring $\triangle^{(n)}\{a,b,c\}$. Thus we find $\partial_{(a,b,c)}(\beta) = \alpha\beta + \beta\alpha = \alpha\beta + \beta$. Similarly $\partial_{(a,b,c)}(\gamma) = \alpha\gamma + \gamma$ and $\partial_{(a,b,c)}(\beta\gamma) = \alpha\beta\gamma + \beta\gamma$. Hence
$$\partial_{(a,b,c)}(\beta)\gamma + \beta\partial_{(a,b,c)}(\gamma) = (\alpha\beta + \beta)\gamma + \beta(\alpha\gamma + \gamma) = \alpha\beta\gamma + \beta\gamma + \beta\alpha\gamma + \beta\gamma =$$ $$= \alpha\beta\gamma + \beta\gamma + \beta\gamma + \beta\gamma = \alpha\beta\gamma + \beta\gamma = \partial_{(a,b,c)}(\beta\gamma).$$
So, we prove
\vspace{3mm}

\textbf{Proposition 1.} \textsl{The map $\partial_{(a,b,c)}$ is a derivation in the whole semiring $\triangle^{(n)}\{a,b,c\}$.}

\vspace{5mm}

\emph{Case 2.} Let $\alpha \in (a,a,a)$. Then $\partial_{(a,a,a)}(\beta) = \alpha\beta + \overline{a} = \alpha\beta$. Similarly $\partial_{(a,a,a)}(\gamma) = \alpha\gamma$ and $\partial_{(a,a,a)}(\beta\gamma) = \alpha\beta\gamma$. Now we find
$\partial_{(a,a,a)}(\beta)\gamma + \beta\partial_{(a,a,a)}(\gamma) = $ $$= \alpha\beta\gamma + \beta\alpha\gamma = \alpha\beta\gamma + \overline{a}\gamma = (\alpha\beta + \overline{a})\gamma = \alpha\beta\gamma = \partial_{(a,a,a)}(\beta\gamma).$$
 Hence we prove
\vspace{3mm}

\textbf{Proposition 2.} \textsl{The map $\partial_{(a,a,a)}$ is a derivation in the whole semiring $\triangle^{(n)}\{a,b,c\}$.}

\vspace{5mm}

\emph{Case 3.}  Let $\alpha \in (c,c,c)$. Then  $\partial_{(c,c,c)}(\beta) = \alpha\beta + \beta\alpha = \alpha\beta + \overline{c} = \overline{c}$. Similarly $\partial_{(c,c,c)}(\gamma) = \overline{c}$ and $\partial_{(c,c,c)}(\beta\gamma) = \overline{c}$. Now
$$\partial_{(c,c,c)}(\beta)\gamma + \beta\partial_{(c,c,c)}(\gamma) = \overline{c}\gamma + \beta\overline{c} = \overline{c}\gamma + \overline{c} = \overline{c} = \partial_{(c,c,c)}(\beta\gamma).$$
 Hence we prove
\vspace{3mm}

\textbf{Proposition 3.} \textsl{The map $\partial_{(c,c,c)}$ is a derivation in the whole semiring $\triangle^{(n)}\{a,b,c\}$.}

\vspace{5mm}

\emph{Case 4.}  Let $\alpha \in (b,b,b)$. For arbitrary endomorphism $\beta$ it follows $\beta\alpha = \overline{b}$ and then $\partial_{(b,b,b)}(\beta) = \alpha\beta + \overline{b}$.
\vspace{2mm}

\emph{Case 4.1.} Let $\beta(b) \leq b$ and $\gamma(b) \leq b$. Now, it follows $\alpha\beta \leq \overline{b}$, $\alpha\gamma \leq \overline{b}$ and then $\partial_{(b,b,b)}(\beta) =  \overline{b}$ and $\partial_{(b,b,b)}(\gamma) =  \overline{b}$. Since $\beta\gamma(b) = \gamma(\beta(b)) \leq \gamma(b) \leq b$, it follows $\partial_{(b,b,b)}(\beta\gamma) =  \overline{b}$. We obtain $\partial_{(b,b,b)}(\beta)\gamma + \beta\partial_{(b,b,b)}(\gamma) = \overline{b}\gamma + \beta\overline{b} = \overline{b}\gamma + \overline{b}$. But $\overline{b}\gamma \leq \overline{b}$ (since $\gamma(b) \leq b$), hence $\partial_{(b,b,b)}(\beta)\gamma + \beta\partial_{(b,b,b)}(\gamma) = \overline{b} = \partial_{(b,b,b)}(\beta\gamma)$.
\vspace{2mm}

\emph{Case 4.2.} Let $\beta(b) \leq b$ and $\gamma \in (c,c,c)$. As in the previous case $\partial_{(b,b,b)}(\beta) =  \overline{b}$. Now, we have $\partial_{(b,b,b)}(\gamma) = \alpha\gamma + \gamma\alpha =   \overline{c} + \overline{b} = \overline{c}$. We obtain $\partial_{(b,b,b)}(\beta\gamma) = \alpha\beta\gamma + \beta\gamma\alpha = \overline{c} + \overline{b} = \overline{c}$ and then $\partial_{(b,b,b)}(\beta)\gamma + \beta\partial_{(b,b,b)}(\gamma) = \overline{b}\gamma + \beta\overline{c} = \overline{c} = \partial_{(b,b,b)}(\beta\gamma)$.
\vspace{2mm}

\emph{Case 4.3.} Let  $\beta \in (c,c,c)$ and $\gamma(b) \leq b$.  As in the previous case $\partial_{(b,b,b)}(\beta) = \overline{c}$ and $\partial_{(b,b,b)}(\gamma) =  \overline{b}$. Now, it follows $\partial_{(b,b,b)}(\beta\gamma) = \alpha\beta\gamma + \beta\gamma\alpha = \overline{c}\gamma + \beta\overline{b}$ and then $\partial_{(b,b,b)}(\beta)\gamma + \beta\partial_{(b,b,b)}(\gamma) = \overline{c}\gamma + \beta\overline{b} = \partial_{(b,b,b)}(\beta\gamma)$.
\vspace{2mm}

\emph{Case 4.4.} Let  $\beta \in (c,c,c)$ and $\gamma \in (c,c,c)$. Obviously $\partial_{(b,b,b)}(\beta) = \overline{c}$, $\partial_{(b,b,b)}(\gamma) =  \overline{c}$, $\partial_{(b,b,b)}(\beta\gamma) = \overline{c}$ and $\partial_{(b,b,b)}(\beta)\gamma + \beta\partial_{(b,b,b)}(\gamma) = \overline{c} = \partial_{(b,b,b)}(\beta\gamma)$.

\vspace{3mm}

Suppose that $\beta(b) = a$, $\gamma(a) \leq b$ and $\gamma(b) = c$. Then $\alpha\beta = \overline{a}$, $\beta\alpha = \overline{b}$ and $\partial_{(b,b,b)}(\beta) = \overline{b}$. We find $\gamma\alpha = \overline{b}$, $\alpha\gamma = \overline{c}$ and $\partial_{(b,b,b)}(\gamma) = \overline{c}$. Now obtain $\partial_{(b,b,b)}(\beta\gamma) = \alpha\beta\gamma + \beta\gamma\alpha = \overline{a}\gamma + \beta\overline{b} = \overline{a}\gamma + \overline{b} = \overline{b}$. Now, it follows
$$\partial_{(b,b,b)}(\beta)\gamma + \beta\partial_{(b,b,b)}(\gamma) = \overline{b}\gamma + \beta\overline{c} = \overline{c} > \partial_{(b,b,b)}(\beta\gamma).$$

From the last inequality follows that the endomorphisms of types $(a,c.c)$ and $(b,c,c)$ does not belong to the subsemiring $\mathcal{D}_{(b,b,b)}$ of $\triangle^{(n)}\{a,b,c\}$, closed under the derivation $\partial_{(b,b,b)}$.

 Hence we prove
 \vspace{3mm}

\textbf{Proposition 4.} \textsl{The map $\partial_{(b,b,b)}$ is a derivation. Maximal subsemiring of  $\triangle^{(n)}\{a,b,c\}$ closed under this derivation is $\mathcal{D}_{(b,b,b)}$ containing all endomorphisms except the endomorphisms of types $(a,c.c)$ and $(b,c,c)$.}
\vspace{5mm}

\emph{Case 5.}  Let $\alpha \in (a,b,b)$.
\vspace{2mm}

\emph{Case 5.1.} Let $\beta(b) \leq b$ and $\gamma(b) \leq b$. Since
$$\begin{array}{c} \alpha\beta(a) = \beta(\alpha(a)) = \beta(a) = \alpha(\beta(a)) = \beta\alpha(a),\\
\alpha\beta(b) = \beta(\alpha(b)) = \beta(b) = \alpha(\beta(b)) = \beta\alpha(b),\\
\alpha\beta(c) = \beta(\alpha(c)) = \beta(b) = \alpha(\beta(b))  \geq \alpha(\beta(c)) = \beta\alpha(c),
\end{array}$$
it follows $\alpha\beta \leq \beta\alpha$, so $\partial_{(a,b,b)}(\beta) = \beta\alpha$. Similarly $\partial_{(a,b,b)}(\gamma) = \gamma\alpha$. Since $\beta\gamma(b) = \gamma(\beta(b) \leq \gamma(b) \leq b$ we have $\partial_{(a,b,b)}(\beta\gamma) = \beta\gamma\alpha$. Now, it follows
$$\partial_{(a,b,b)}(\beta)\gamma + \beta\partial_{(a,b,b)}(\gamma) = \beta\alpha\gamma + \beta\gamma\alpha = \beta(\alpha\gamma + \gamma\alpha) = \beta\gamma\alpha = \partial_{(a,b,b)}(\beta\gamma).$$

\emph{Case 5.2.} Let $\beta(b) \leq b$ and $\gamma \in (c,c,c)$. As in the previous case follows $\partial_{(a,b,b)}(\beta) = \beta\alpha$. We obtain $\partial_{(a,b,b)}(\gamma) = \alpha\gamma + \gamma\alpha = \overline{c} + \gamma\alpha = \overline{c}$. So, we find $\partial_{(a,b,b)}(\beta\gamma) = \alpha\beta\gamma + \beta\gamma\alpha = \overline{c} + \beta\gamma\alpha = \overline{c}$ and
$$\partial_{(a,b,b)}(\beta)\gamma + \beta\partial_{(a,b,b)}(\gamma) = \beta\alpha\gamma + \beta\gamma\alpha = \beta(\alpha\gamma + \gamma\alpha) = \beta\overline{c} = \overline{c} = \partial_{(a,b,b)}(\beta\gamma).$$

\emph{Case 5.3.} Let $\beta \in (c,c,c)$ and $\beta(b) \leq b$. As in the previous case $\alpha\beta = \overline{c}$, $\partial_{(a,b,b)}(\beta) = \overline{c}$ and $\partial_{(a,b,b)}(\gamma) = \gamma\alpha$. We obtain $\partial_{(a,b,b)}(\beta\gamma) = \alpha\beta\gamma + \beta\gamma\alpha = \overline{c}\gamma + \beta\gamma\alpha$ and then
$$\partial_{(a,b,b)}(\beta)\gamma + \beta\partial_{(a,b,b)}(\gamma) = \overline{c}\gamma + \beta\gamma\alpha = \partial_{(a,b,b)}(\beta\gamma).$$

\emph{Case 5.4.} Let  $\beta \in (c,c,c)$ and $\gamma \in (c,c,c)$. Obviously $\partial_{(a,b,b)}(\beta) = \overline{c}$, $\partial_{(a,b,b)}(\gamma) =  \overline{c}$, $\partial_{(a,b,b)}(\beta\gamma) = \overline{c}$ and $\partial_{(a,b,b)}(\beta)\gamma + \beta\partial_{(a,b,b)}(\gamma) = \overline{c} = \partial_{(a,b,b)}(\beta\gamma)$.
\vspace{4mm}

Suppose that $\beta(b) = a$, $\beta(c) = b$, $\gamma(a) = b$  and $\gamma(b) = c$. Since
$$\begin{array}{c} \alpha\beta(a) = \beta(\alpha(a)) = \beta(a) = a = \alpha(\beta(a)) = \beta\alpha(a),\\
\alpha\beta(b) = \beta(\alpha(b)) = \beta(b) = a = \alpha(\beta(b)) = \beta\alpha(b),\\
\alpha\beta(c) = \beta(\alpha(c)) = \beta(b) = a < b = \alpha(\beta(c)) = \beta\alpha(c),
\end{array}$$
it follows $\partial_{(a,b,b)}(\beta) = \beta\alpha$. Since
$$\begin{array}{c} \alpha\gamma(a) = \gamma(\alpha(a)) = \gamma(a) = b = \alpha(\gamma(a)) = \gamma\alpha(a),\\
\alpha\gamma(b) = \gamma(\alpha(b)) = \gamma(b) = c > b = \alpha(\gamma(b)) = \gamma\alpha(b),\\
\alpha\gamma(c) = \gamma(\alpha(c)) = \gamma(b) = c > b = \alpha(\gamma(c)) = \gamma\alpha(c),
\end{array}$$
it follows $\partial_{(a,b,b)}(\gamma) = \alpha\gamma$. We find $\partial_{(a,b,b)}(\beta\gamma) = \alpha\beta\gamma + \beta\gamma\alpha < \beta\alpha\gamma + \beta\alpha\gamma = \beta\alpha\gamma$ and then
$$\partial_{(a,b,b)}(\beta)\gamma + \beta\partial_{(a,b,b)}(\gamma) = \beta\alpha\gamma + \beta\alpha\gamma = \beta\alpha\gamma > \partial_{(a,b,b)}(\beta\gamma).$$

From the last inequality follows that the endomorphisms of type  $(b,c,c)$ does not belong to the subsemiring $\mathcal{D}_{(a,b,b)}$ of $\triangle^{(n)}\{a,b,c\}$, closed under the derivation $\partial_{(a,b,b)}$. If  assume that endomorphisms of type $(a,c,c)$ belongs to $\mathcal{D}_{(a,b,b)}$ we have reached a contradiction, because sum of an endomorphism of type $(a,c,c)$ and $\overline{b}$ is an endomorphism of type $(b,c,c)$. So, $\mathcal{D}_{(a,b,b)} = \mathcal{D}_{(b,b,b)}$. Thus we prove
\vspace{2mm}

\textbf{Proposition 5.} \textsl{The map $\partial_{(a,b,b)}$ is a derivation. Maximal subsemiring of  $\triangle^{(n)}\{a,b,c\}$ closed under this derivation is $\mathcal{D}_{(a,b,b)}$ containing all endomorphisms except the endomorphisms of types $(a,c.c)$ and $(b,c,c)$.}
\vspace{5mm}

\emph{Case 6.}  Let $\alpha \in (a,a,c)$.
\vspace{2mm}

\emph{Case 6.1.} Let $\beta(b) \leq b$ and $\gamma(b) \leq b$. Since
$$\begin{array}{c} \alpha\beta(a) = \beta(\alpha(a)) = \beta(a) \geq a  = \alpha(\beta(a)) = \beta\alpha(a),\\
\alpha\beta(b) = \beta(\alpha(b)) = \beta(a) \geq a = \alpha(\beta(b)) = \beta\alpha(b),\\
\alpha\beta(c) = \beta(\alpha(c)) = \beta(c)  = \alpha(\beta(c)) = \beta\alpha(c),
\end{array}$$
it follows $\alpha\beta \geq \beta\alpha$ and $\partial_{(a,a,c)}(\beta) = \alpha\beta$. Similarly $\partial_{(a,a,c)}(\gamma) = \alpha\gamma$. Since  $\beta\gamma(b) = \gamma(\beta(b) \leq \gamma(b) \leq b$ we have $\partial_{(a,a,c)}(\beta\gamma) = \alpha\beta\gamma$. So, we have
$$\partial_{(a,a,c)}(\beta)\gamma  + \beta\partial_{(a,a,c)}(\gamma) = \alpha\beta\gamma + \beta\alpha\gamma = (\alpha\beta + \beta\alpha)\gamma = \alpha\beta\gamma = \partial_{(a,a,c)}(\beta\gamma).$$

\emph{Case 6.2.} Let $\beta(b) \leq b$ and $\gamma \in (c,c,c)$. As in the previous case follows $\partial_{(a,a,c)}(\beta) = \alpha\beta$. It follows $\partial_{(a,a,c)}(\gamma) = \alpha\gamma + \gamma\alpha  = \overline{c} + \overline{c}  = \overline{c}$. So, we find $\partial_{(a,a,c)}(\beta\gamma) = \alpha\beta\gamma + \beta\gamma\alpha = \overline{c} + \beta\gamma\alpha = \overline{c}$ and then
$$\partial_{(a,a,c)}(\beta)\gamma + \beta\partial_{(a,a,c)}(\gamma) = \alpha\beta\gamma + \beta\overline{c} = \overline{c} + \overline{c} =  \overline{c} = \partial_{(a,a,c)}(\beta\gamma).$$

\emph{Case 6.3.} Let $\beta \in (c,c,c)$ and $\gamma(b) \leq b$. As in the previous case $\alpha\beta = \overline{c}$ and $\beta\alpha = \overline{c}$, so, $\partial_{(a,a,c)}(\beta) = \overline{c}$ and also $\partial_{(a,a,c)}(\gamma) = \alpha\gamma$. Now, we obtain $\partial_{(a,a,c)}(\beta\gamma) = \alpha\beta\gamma + \beta\gamma\alpha = \overline{c}\gamma + \beta\gamma\alpha$. But $\beta\gamma\alpha \leq \beta\alpha\gamma = \overline{c}\gamma$, hence, $\partial_{(a,a,c)}(\beta\gamma) = \overline{c}\gamma$. Then
$$\partial_{(a,a,c)}(\beta)\gamma + \beta\partial_{(a,a,c)}(\gamma) = \overline{c}\gamma + \beta\alpha\gamma = \overline{c}\gamma + \overline{c}\gamma = \overline{c}\gamma = \partial_{(a,a,c)}(\beta\gamma).$$

\emph{Case 6.4.} Let  $\beta \in (c,c,c)$ and $\gamma \in (c,c,c)$. As in case 4 and case 5 we find $\partial_{(a,a,c)}(\beta)\gamma + \beta\partial_{(a,a,c)}(\gamma) = \partial_{(a,a,c)}(\beta\gamma)$.
\vspace{5mm}

Suppose that $\beta(a) = a$, $\beta(b) = c$  and $\gamma(b) \leq b$. Since
$$\begin{array}{c} \alpha\beta(a) = \beta(\alpha(a)) = \beta(a) = a = \alpha(\beta(a)) = \beta\alpha(a),\\
\alpha\beta(b) = \beta(\alpha(b)) = \beta(a) = a < c = \alpha(\beta(b)) = \beta\alpha(b),\\
\alpha\beta(c) = \beta(\alpha(c)) = \beta(c) = c = \alpha(\beta(c)) = \beta\alpha(c),
\end{array}$$
it follows $\alpha\beta < \beta\alpha$ and $\partial_{(a,a,c)}(\beta) = \beta\alpha$. As in case 6.1
$\partial_{(a,a,c)}(\gamma) = \alpha\gamma$. We calculate $\partial_{(a,a,c)}(\beta)\gamma + \beta\partial_{(a,a,c)}(\gamma) = \beta\alpha\gamma + \beta\alpha\gamma = \beta\alpha\gamma$. Now, it follows
$$\partial_{(a,a,c)}(\beta\gamma) = \alpha\beta\gamma + \beta\gamma\alpha > \beta\alpha\gamma + \beta\gamma\alpha = $$ $$= \beta(\alpha\gamma + \gamma\alpha) =  \beta\alpha\gamma = \partial_{(a,a,c)}(\beta)\gamma + \beta\partial_{(a,a,c)}(\gamma).$$

From the last inequality follows that the endomorphisms of type  $(a,c,c)$ does not belong to the subsemiring $\mathcal{D}_{(a,a,c)}$ of $\triangle^{(n)}\{a,b,c\}$, closed under the derivation $\partial_{(a,a,c)}$. If  assume that endomorphisms of type $(b,c,c)$ belongs to $\mathcal{D}_{(a,a,c)}$ we find a contradiction, because product of an endomorphism of type $(b,c,c)$ and an endomorphism of type $(a,a,c)$  is an endomorphism of type $(a,c,c)$. So, $\mathcal{D}_{(a,a,c)} = \mathcal{D}_{(b,b,b)}$. Thus we prove
\vspace{3mm}

\textbf{Proposition 6.} \textsl{The map $\partial_{(a,a,c)}$ is a derivation. Maximal subsemiring of  $\triangle^{(n)}\{a,b,c\}$ closed under this derivation is $\mathcal{D}_{(a,a,c)}$ containing all endomorphisms except the endomorphisms of types $(a,c.c)$ and $(b,c,c)$.}
\vspace{5mm}

\emph{Case 7.}  Let $\alpha \in (b,b,c)$.
\vspace{1mm}

\emph{Case 7.1.} Let $\beta(b) \leq b$ and $\gamma(b) \leq b$. Since
$$\begin{array}{c} \alpha\beta(a) = \beta(\alpha(a)) = \beta(b) \leq b  = \alpha(\beta(a)) = \beta\alpha(a),\\
\alpha\beta(b) = \beta(\alpha(b)) = \beta(b) \leq b = \alpha(\beta(b)) = \beta\alpha(b),\\
\alpha\beta(c) = \beta(\alpha(c)) = \beta(c)  \leq \alpha(\beta(c)) = \beta\alpha(c),
\end{array}$$
it follows $\alpha\beta \leq \beta\alpha$ and $\partial_{(b,b,c)}(\beta) = \beta\alpha$. Similarly $\partial_{(b,b,c)}(\gamma) = \gamma\alpha$. Since  $\beta\gamma(b) = \gamma(\beta(b) \leq \gamma(b) \leq b$, it follows $\partial_{(b,b,c)}(\beta\gamma) = \beta\gamma\alpha$. So, we obtain
$$\partial_{(b,b,c)}(\beta)\gamma  + \beta\partial_{(b,b,c)}(\gamma) = \beta\alpha\gamma + \beta\gamma\alpha = \beta(\alpha\gamma + \gamma\alpha) = \beta\gamma\alpha = \partial_{(b,b,c)}(\beta\gamma).$$

\emph{Case 7.2.} Let $\beta(b) \leq b$ and $\gamma \in (c,c,c)$. As in the previous case follows $\partial_{(b,b,c)}(\beta) = \beta\alpha$. We find $\partial_{(b,b,c)}(\gamma) = \alpha\gamma + \gamma\alpha  = \overline{c} + \overline{c}  = \overline{c}$. Then we find $\partial_{(b,b,c)}(\beta\gamma) = \alpha\beta\gamma + \beta\gamma\alpha = \overline{c} + \overline{c}\alpha = \overline{c}$. So, it follows
$$\partial_{(b,b,c)}(\beta)\gamma + \beta\partial_{(b,b,c)}(\gamma) = \beta\alpha\gamma + \beta\overline{c} = \overline{c} + \overline{c} =  \overline{c} = \partial_{(b,b,c)}(\beta\gamma).$$

\emph{Case 7.3.} Let $\beta \in (c,c,c)$ and $\gamma(b) \leq b$. As in the previous case $\alpha\beta = \overline{c}$ and $\beta\alpha = \overline{c}$, then $\partial_{(b,b,c)}(\beta) = \overline{c}$ and also $\partial_{(b,b,c)}(\gamma) = \gamma\alpha$. Now, we obtain $\partial_{(b,b,c)}(\beta\gamma) = \alpha\beta\gamma + \beta\gamma\alpha = \overline{c}\gamma + \beta\gamma\alpha$. Then
$$\partial_{(b,b,c)}(\beta)\gamma + \beta\partial_{(b,b,c)}(\gamma) = \overline{c}\gamma + \beta\gamma\alpha =  \partial_{(b,b,c)}(\beta\gamma).$$

\emph{Case 7.4.} Let  $\beta \in (c,c,c)$ and $\gamma \in (c,c,c)$. As in the cases 4, 5 and 6 we find $\partial_{(b,b,c)}(\beta)\gamma + \beta\partial_{(b,b,c)}(\gamma) = \partial_{(b,b,c)}(\beta\gamma)$.
\vspace{5mm}

Suppose that $\beta \in (a,a,b)$  and $\gamma \in (a,c,c)$. Since
$$\begin{array}{c} \alpha\beta(a) = \beta(\alpha(a)) = \beta(b) = a < b = \alpha(\beta(a)) = \beta\alpha(a),\\
\alpha\beta(b) = \beta(\alpha(b)) = \beta(b) = a < b = \alpha(\beta(b)) = \beta\alpha(b),\\
\alpha\beta(c) = \beta(\alpha(c)) = \beta(c) = b = \alpha(\beta(c)) = \beta\alpha(c),
\end{array}$$
it follows $\alpha\beta < \beta\alpha$ and $\partial_{(b,b,c)}(\beta) = \beta\alpha$. Since
$$\begin{array}{c} \alpha\gamma(a) = \gamma(\alpha(a)) = \gamma(b) = c > b = \alpha(\gamma(a)) = \gamma\alpha(a),\\
\alpha\gamma(b) = \gamma(\alpha(b)) = \gamma(b) = c  = \alpha(\gamma(b)) = \gamma\alpha(b),\\
\alpha\gamma(c) = \gamma(\alpha(c)) = c > b = \alpha(\gamma(c)) = \gamma\alpha(c),
\end{array}$$
it follows $\alpha\gamma > \gamma\alpha$ and $\partial_{(b,b,c)}(\gamma) = \alpha\gamma$

 Now, it follows $\partial_{(b,b,c)}(\beta)\gamma + \beta\partial_{(b,b,c)}(\gamma) = \beta\alpha\gamma + \beta\alpha\gamma = \beta\alpha\gamma$. Then
 $$\partial_{(b,b,c)}(\beta\gamma) = \alpha\beta\gamma + \beta\gamma\alpha < \beta\alpha\gamma + \beta\gamma\alpha = $$ $$= \beta(\alpha\gamma + \gamma\alpha) =  \beta\alpha\gamma = \partial_{(b,b,c)}(\beta)\gamma + \beta\partial_{(b,b,c)}(\gamma).$$

 From the last inequality follows that the endomorphisms of type  $(a,c,c)$ does not belong to the subsemiring $\mathcal{D}_{(b,b,c)}$ of $\triangle^{(n)}\{a,b,c\}$, closed under the derivation $\partial_{(b,b,c)}$. If  assume that endomorphisms of type $(b,c,c)$ belongs to $\mathcal{D}_{(b,b,c)}$ we find a contradiction, because product of an endomorphism of type $(b,c,c)$ and an endomorphism of type $(a,a,c)$  is an endomorphism of type $(a,c,c)$. So, $\mathcal{D}_{(b,b,c)} = \mathcal{D}_{(b,b,b)}$. Thus we prove
\vspace{3mm}

\textbf{Proposition 7.} \textsl{The map $\partial_{(b,b,c)}$ is a derivation. Maximal subsemiring of  $\triangle^{(n)}\{a,b,c\}$ closed under this derivation is $\mathcal{D}_{(b,b,c)}$ containing all endomorphisms except the endomorphisms of types $(a,c,c)$ and $(b,c,c)$.}
\vspace{5mm}

\emph{Case 8.}  Let $\alpha \in (a,c,c)$.
c

\emph{Case 8.1.} Let $\beta(b) \geq b$ and $\gamma(b) \geq b$. Since
$$\begin{array}{c} \alpha\beta(a) = \beta(\alpha(a)) = \beta(a) \leq \alpha(\beta(a)) = \beta\alpha(a),\\
\alpha\beta(b) = \beta(\alpha(b)) = \beta(c) \leq c = \alpha(\beta(b)) = \beta\alpha(b),\\
\alpha\beta(c) = \beta(\alpha(c)) = \beta(c)  \leq c = \alpha(\beta(c)) = \beta\alpha(c),
\end{array}$$
it follows $\alpha\beta \leq \beta\alpha$ and $\partial_{(a,c,c)}(\beta) = \beta\alpha$. Similarly $\partial_{(a,c,c)}(\gamma) = \gamma\alpha$. Since  $\beta\gamma(b) = \gamma(\beta(b) \geq \gamma(b) \geq b$, it follows $\partial_{(a,c,c)}(\beta\gamma) = \beta\gamma\alpha$. So, we obtain
$$\partial_{(a,c,c)}(\beta)\gamma  + \beta\partial_{(a,c,c)}(\gamma) = \beta\alpha\gamma + \beta\gamma\alpha = \beta(\alpha\gamma + \gamma\alpha) = \beta\gamma\alpha = \partial_{(a,c,c)}(\beta\gamma).$$

\emph{Case 8.2.} Let $\beta(b) \geq b$ and $\gamma \in (a,a,a)$. As in the previous case follows $\partial_{(a,c,c)}(\beta) = \beta\alpha$. We obtain $\partial_{(a,c,c)}(\gamma) = \alpha\gamma + \gamma\alpha  = \overline{a} + \overline{a}  = \overline{a}$. Then we find $\partial_{(a,c,c)}(\beta\gamma) = \alpha\beta\gamma + \beta\gamma\alpha = \overline{a} + \beta\overline{a} = \overline{a}$. So, it follows
$$\partial_{(a,c,c)}(\beta)\gamma + \beta\partial_{(a,c,c)}(\gamma) = \beta\alpha\gamma + \beta\overline{a} = \overline{a} + \overline{a} =  \overline{a} = \partial_{(a,c,c)}(\beta\gamma).$$

\emph{Case 8.3.} Let $\beta \in (a,a,a)$ and $\gamma(b) \geq b$. As in the previous case $\alpha\beta = \overline{a}$ and $\beta\alpha = \overline{a}$, then $\partial_{(a,c,c)}(\beta) = \overline{a}$ and also $\partial_{(a,c,c)}(\gamma) = \gamma\alpha$. Now, we find $\partial_{(a,c,c)}(\beta\gamma) = \alpha\beta\gamma + \beta\gamma\alpha = \overline{a}\gamma + \beta\gamma\alpha$. Then
$$\partial_{(a,c,c)}(\beta)\gamma + \beta\partial_{(a,c,c)}(\gamma) = \overline{a}\gamma + \beta\gamma\alpha =  \partial_{(a,c,c)}(\beta\gamma).$$

\emph{Case 8.4.} Let  $\beta \in (a,a,a)$ and $\gamma \in (a,a,a)$. As in the previous case $\partial_{(a,c,c)}(\beta) = \overline{a}$, $\partial_{(a,c,c)}(\gamma) = \overline{a}$ and $\partial_{(a,c,c)}(\beta\gamma) = \overline{a}$. Then we have $\partial_{(a,c,c)}(\beta\gamma) = \partial_{(a,c,c)}(\beta)\gamma + \beta\partial_{(a,c,c)}(\gamma)$.
\vspace{6mm}

Suppose that $\beta \in (a,b,b)$  and $\gamma \in (a,a,c)$. Since
$$\begin{array}{c} \alpha\beta(a) = \beta(\alpha(a)) = \beta(a) = a  = \alpha(\beta(a)) = \beta\alpha(a),\\
\alpha\beta(b) = \beta(\alpha(b)) = \beta(c) = b < c = \alpha(\beta(b)) = \beta\alpha(b),\\
\alpha\beta(c) = \beta(\alpha(c)) = \beta(c) = b < c = \alpha(\beta(c)) = \beta\alpha(c),
\end{array}$$
it follows $\alpha\beta < \beta\alpha$ and $\partial_{(a,c,c)}(\beta) = \beta\alpha$. Since
$$\begin{array}{c} \alpha\gamma(a) = \gamma(\alpha(a)) = \gamma(a) = a = \alpha(\gamma(a)) = \gamma\alpha(a),\\
\alpha\gamma(b) = \gamma(\alpha(b)) = \gamma(c) = c > a = \alpha(\gamma(b)) = \gamma\alpha(b),\\
\alpha\gamma(c) = \gamma(\alpha(c)) = c  = \alpha(\gamma(c)) = \gamma\alpha(c),
\end{array}$$
it follows $\alpha\gamma > \gamma\alpha$ and $\partial_{(a,c,c)}(\gamma) = \alpha\gamma$
\vspace{3mm}

 Now, it follows $\partial_{(a,c,c)}(\beta)\gamma + \beta\partial_{(a,c,c)}(\gamma) = \beta\alpha\gamma + \beta\alpha\gamma = \beta\alpha\gamma$. Then
 $$\partial_{(a,c,c)}(\beta\gamma) = \alpha\beta\gamma + \beta\gamma\alpha < \beta\alpha\gamma + \beta\gamma\alpha = $$ $$= \beta(\alpha\gamma + \gamma\alpha) =  \beta\alpha\gamma = \partial_{(a,c,c)}(\beta)\gamma + \beta\partial_{(b,b,c)}(\gamma).$$

 From the last inequality follows that the endomorphisms of type  $(a,a,c)$ does not belong to the subsemiring $\mathcal{D}_{(a,c,c)}$ of $\triangle^{(n)}\{a,b,c\}$, closed under the derivation $\partial_{(a,c,c)}$.
 \vspace{3mm}

  If  assume that endomorphisms of type $(a,a,b)$ belongs to $\mathcal{D}_{(a,c,c)}$ we find a contradiction, because product of an endomorphism of type $(a,a,b)$ and an endomorphism of type $(a,c,c)$  is an endomorphism of type $(a,a,c)$.
\vspace{3mm}

 So, the semiring $\mathcal{D}_{(a,c,c)}$ does not contain endomorphisms of types $(a,a,b)$ and $(a,a,c)$. Thus we prove
\vspace{4mm}

\textbf{Proposition 8.} \textsl{The map $\partial_{(a,c,c)}$ is a derivation. Maximal subsemiring of  $\triangle^{(n)}\{a,b,c\}$ closed under this derivation is $\mathcal{D}_{(a,c,c)}$ containing all endomorphisms except the endomorphisms of types $(a,a,b)$ and $(a,a,c)$.}

\emph{Case 9.}  Let $\alpha \in (a,a,b)$.
\vspace{2mm}

\emph{Case 9.1.} Let $\beta(b) = a $ and $\gamma(b) = a$. Then $\alpha\beta = \overline{a}$ and $\partial_{(a,a,b)}(\beta) = \beta\alpha$. Similarly $\alpha\gamma = \overline{a}$ and $\partial_{(a,a,b)}(\gamma) = \gamma\alpha$. Since $\beta\gamma(b) = \gamma(\beta(b)) = \gamma(b) = b$, it follows $\partial_{(a,a,b)}(\beta\gamma) = \beta\gamma\alpha$. Hence
$$\partial_{(a,a,b)}(\beta)\gamma + \beta\partial_{(a,a,b)}(\gamma) = \beta\alpha\gamma + \beta\gamma\alpha = \beta(\alpha\gamma + \gamma\alpha)  = \beta\gamma\alpha = \partial_{(a,a,b)}(\beta\gamma).$$

\emph{Case 9.2.} Let $\beta(b) = b$ and $\gamma(b) = b$. Since
$$\begin{array}{c} \alpha\beta(a) = \beta(\alpha(a)) = \beta(a) \geq \alpha(\beta(a)) = \beta\alpha(a),\\
\alpha\beta(b) = \beta(\alpha(b)) = \beta(a) \geq a = \alpha(b) = \alpha(\beta(b)) = \beta\alpha(b),\\
\alpha\beta(c) = \beta(\alpha(c)) = \beta(b) = b = \alpha(c) = \alpha(\beta(c)) = \beta\alpha(c),
\end{array}$$
it follows $\alpha\beta \geq \beta\alpha$ and $\partial_{(a,a,b)}(\beta) = \alpha\beta$. Similarly $\partial_{(a,a,b)}(\gamma) = \alpha\gamma$. Since $\beta\gamma(b) = b$, it follows $\partial_{(a,a,b)}(\beta\gamma) = \alpha\beta\gamma$. Then we obtain
$$\partial_{(a,a,b)}(\beta)\gamma + \beta\partial_{(a,a,b)}(\gamma) = \alpha\beta\gamma + \beta\alpha\gamma = (\alpha\beta + \beta\alpha)\gamma = \alpha\beta\gamma =  \partial_{(a,a,b)}(\beta\gamma).$$

\emph{Case 9.3.} Let $\beta(b) = a$ and $\gamma(b) = b$. As in case 9.1, it follows $\alpha\beta = \overline{a}$ and $\partial_{(a,a,b)}(\beta) = \beta\alpha$. As in case 9.2 we obtain $\partial_{(a,a,b)}(\gamma) = \gamma\alpha$.
We calculate $\partial_{(a,a,b)}(\beta)\gamma + \beta\partial_{(a,a,b)}(\gamma) = \beta\alpha\gamma + \beta\gamma\alpha = \beta(\alpha\gamma + \gamma\alpha) = \beta\gamma\alpha$ and $\partial_{(a,a,b)}(\beta\gamma) = \alpha\beta\gamma + \beta\gamma\alpha = \overline{a}\gamma + \beta\gamma\alpha$. Since
$\beta\gamma\alpha \geq \beta\alpha\gamma \geq \overline{a}\gamma$, it follows  $$\partial_{(a,a,b)}(\beta\gamma)  = \beta\gamma\alpha = \partial_{(a,a,b)}(\beta)\gamma + \beta\partial_{(a,a,b)}(\gamma).$$

\emph{Case 9.4.} Let $\beta(b) = b$ and $\gamma(b) = a$. As in the previous case  $\partial_{(a,a,b)}(\beta) = \beta\alpha$, $\partial_{(a,a,b)}(\gamma) = \gamma\alpha$ and $\partial_{(a,a,b)}(\beta)\gamma + \beta\partial_{(a,a,b)}(\gamma) = \beta\alpha\gamma + \beta\gamma\alpha = \beta(\alpha\gamma + \gamma\alpha) = \beta\gamma\alpha$. Since
$\alpha\beta\gamma \leq \beta\alpha\gamma$, it follows  $$\partial_{(a,a,b)}(\beta\gamma) = \alpha\beta\gamma + \beta\gamma\alpha  = \beta\gamma\alpha = \partial_{(a,a,b)}(\beta)\gamma + \beta\partial_{(a,a,b)}(\gamma).$$

\emph{Case 9.5.} Let $\beta(b) \leq b$ and $\gamma \in (c,c,c)$. We obtain
$\partial_{(a,a,b)}(\gamma) = \alpha\gamma + \gamma\alpha  = \overline{c} + \gamma\alpha  = \overline{c}$ and $\partial_{(a,a,b)}(\beta\gamma) = \alpha\beta\gamma + \beta\gamma\alpha = \overline{c} + \beta\gamma\alpha = \overline{c}$. So, it follows
$$\partial_{(a,a,b)}(\beta)\gamma + \beta\partial_{(a,a,b)}(\gamma) = (\alpha\beta +\beta\alpha)\gamma + \beta\overline{c} = (\alpha\beta +\beta\alpha)\gamma + \overline{c} =  \overline{c} = \partial_{(b,b,c)}(\beta\gamma).$$

\emph{Case 9.6.} Let $\beta \in (c,c,c)$ and $\gamma(b) \leq b$. As in the previous case $\alpha\beta = \overline{c}$ and  $\partial_{(a,a,b)}(\beta) = \overline{c}$. Now, we obtain $\partial_{(a,a,b)}(\beta\gamma) = \alpha\beta\gamma + \beta\gamma\alpha = \overline{c}\gamma + \beta\gamma\alpha$.
Then
$$\partial_{(a,a,b)}(\beta)\gamma + \beta\partial_{(a,a,b)}(\gamma) = \overline{c}\gamma + \beta(\alpha\gamma + \gamma\alpha) = \overline{c}\gamma + \beta\alpha\gamma + \beta\gamma\alpha = $$ $$= (\overline{c} + \beta\alpha)\gamma + \beta\gamma\alpha = \overline{c}\gamma + \beta\gamma\alpha =  \partial_{(a,a,c)}(\beta\gamma).$$

\emph{Case 9.7.} Let  $\beta \in (c,c,c)$ and $\gamma \in (c,c,c)$. As in the previous cases  we find $\partial_{(a,a,b)}(\beta) = \overline{c}$, $\partial_{(a,a,b)}(\gamma) = \overline{c}$ and $\partial_{(a,a,b)}(\beta\gamma) = \overline{c}$. Then, it follows $$\partial_{(a,a,b)}(\beta)\gamma + \beta\partial_{(a,a,b)}(\gamma) = \overline{c} = \partial_{(a,a,b)}(\beta\gamma).$$
\vspace{2mm}

Suppose that $\beta(b) = a$,  $\gamma(a) = b$ and $\gamma(b) = c$.
As in case 9.1, it follows $\alpha\beta = \overline{a}$ and $\partial_{(a,a,b)}(\beta) = \beta\alpha$. Since
$$\begin{array}{c} \alpha\gamma(a) = \gamma(\alpha(a)) = \gamma(a) = b > a = \gamma(b) = \alpha(\gamma(a)) = \gamma\alpha(a),\\
\alpha\gamma(b) = \gamma(\alpha(b)) = \gamma(a) = b = \alpha(c) = \alpha(\gamma(b)) = \gamma\alpha(b),\\
\alpha\gamma(c) = \gamma(\alpha(c)) = \gamma(b) = c  > b = \alpha(c) = \alpha(\gamma(c)) = \gamma\alpha(c),
\end{array}$$
it follows $\alpha\gamma > \gamma\alpha$ and $\partial_{(a,a,b)}(\gamma) = \alpha\gamma$. We obtain $\partial_{(a,a,b)}(\beta\gamma) = \alpha\beta\gamma + \beta\gamma\alpha = \overline{a}\gamma + \beta\gamma\alpha < \overline{a}\gamma + \beta\alpha\gamma = (\overline{a} + \beta\alpha)\gamma = \beta\alpha\gamma$. Then, it follows
$$\partial_{(a,a,b)}(\beta)\gamma + \beta\partial_{(a,a,b)}(\gamma) = \beta\alpha\gamma + \beta\alpha\gamma = \beta\alpha\gamma > \partial_{(a,a,b)}(\beta\gamma).$$

From the last inequality follows that the endomorphisms of type  $(b,c,c)$ does not belong to the subsemiring $\mathcal{D}_{(a,a,b)}$ of $\triangle^{(n)}\{a,b,c\}$, closed under the derivation $\partial_{(a,a,b)}$. If  assume that endomorphisms of type $(a,c,c)$ belongs to $\mathcal{D}_{(a,a,b)}$ we have reached a contradiction, because sum of an endomorphism of type $(a,c,c)$ and $\overline{b}$ is an endomorphism of type $(b,c,c)$. So, $\mathcal{D}_{(a,a,b)} = \mathcal{D}_{(b,b,b)}$. Thus we prove
\vspace{3mm}

\textbf{Proposition 9.} \textsl{The map $\partial_{(a,a,b)}$ is a derivation. Maximal subsemiring of  $\triangle^{(n)}\{a,b,c\}$ closed under this derivation is $\mathcal{D}_{(a,a,b)}$ containing all endomorphisms except the endomorphisms of types $(a,c.c)$ and $(b,c,c)$.}
\vspace{3mm}

\vspace{5mm}

\emph{Case 10.}  Let $\alpha \in (b,c,c)$.
\vspace{2mm}

\emph{Case 10.1.} Let $\beta(a) \geq b$ and $\gamma(a) \geq b$. Since
$$\begin{array}{c} \alpha\beta(a) = \beta(\alpha(a)) = \beta(b) \leq c = \alpha(b) \leq \alpha(\beta(a)) = \beta\alpha(a),\\
\alpha\beta(b) = \beta(\alpha(b)) = \beta(c) \leq c = \alpha(b) \leq \alpha(\beta(b)) = \beta\alpha(b),\\
\alpha\beta(c) = \beta(\alpha(c)) = \beta(c)  \leq  \alpha(\beta(c)) = \beta\alpha(c),
\end{array}$$
it follows $\alpha\beta \leq \beta\alpha$ and $\partial_{(b,c,c)}(\beta) = \beta\alpha$. Similarly $\partial_{(b,c,c)}(\gamma) = \gamma\alpha$. Since  $\beta\gamma(a) = \gamma(\beta(a) \geq \gamma(b) \geq \gamma(a) \geq b$, it follows $\partial_{(b,c,c)}(\beta\gamma) = \beta\gamma\alpha$. So, we obtain
$$\partial_{(b,c,c)}(\beta)\gamma  + \beta\partial_{(b,c,c)}(\gamma) = \beta\alpha\gamma + \beta\gamma\alpha = \beta(\alpha\gamma + \gamma\alpha) = \beta\gamma\alpha = \partial_{(b,c,c)}(\beta\gamma).$$

\emph{Case 10.2.} Let $\beta(a) = a$, $\beta(b) = b$, $\gamma(a) = a$ and $\gamma(b) = b$. Since
$$\begin{array}{c} \alpha\beta(a) = \beta(\alpha(a)) = \beta(b) = b = \alpha(a) = \alpha(\beta(a)) = \beta\alpha(a),\\
\alpha\beta(b) = \beta(\alpha(b)) = \beta(c) \leq c = \alpha(b) \leq \alpha(\beta(b)) = \beta\alpha(b),\\
\alpha\beta(c) = \beta(\alpha(c)) = \beta(c)  \leq c = \alpha(b) \leq \alpha(\beta(c)) = \beta\alpha(c),
\end{array}$$
it follows $\alpha\beta \leq \beta\alpha$ and $\partial_{(b,c,c)}(\beta) = \beta\alpha$. Similarly $\partial_{(b,c,c)}(\gamma) = \gamma\alpha$. Since  $\beta\gamma$ has the same properties as $\beta$ and $\gamma$,  it follows $\partial_{(b,c,c)}(\beta\gamma) = \beta\gamma\alpha$. As in the previous case we obtain the equality
$$\partial_{(b,c,c)}(\beta)\gamma  + \beta\partial_{(b,c,c)}(\gamma) =  \partial_{(b,c,c)}(\beta\gamma).$$

\emph{Case 10.3.} Let $\beta(a) \geq b$,  $\gamma(a) = a$ and $\gamma(b) = b$. As in case 10.1 we find $\partial_{(b,c,c)}(\beta) = \beta\alpha$. As in case 10.2 we obtain $\partial_{(b,c,c)}(\gamma) = \gamma\alpha$.

 Now we find $\partial_{(b,c,c)}(\beta)\gamma + \beta\partial_{(b,c,c)}(\gamma) = \beta\alpha\gamma + \beta\gamma\alpha = \beta(\alpha\gamma + \gamma\alpha) = \beta\gamma\alpha$
and $\partial_{(a,c,c)}(\beta\gamma) = \alpha\beta\gamma + \beta\gamma\alpha$. Since
$$\beta\gamma\alpha \leq \alpha\beta\gamma + \beta\gamma\alpha \leq \beta\alpha\gamma + \beta\gamma\alpha = \beta(\alpha\gamma + \gamma\alpha) = \beta\gamma\alpha,$$
it follows
$$\partial_{(a,c,c)}(\beta)\gamma + \beta\partial_{(a,c,c)}(\gamma) = \beta\gamma\alpha =  \partial_{(a,c,c)}(\beta\gamma).$$

\emph{Case 10.4.} Let $\beta(a) = a$, $\beta(b) = b$ and $\gamma(a) \geq b$. As in the previous case we find $\partial_{(b,c,c)}(\beta) = \beta\alpha$ and $\partial_{(b,c,c)}(\gamma) = \gamma\alpha$. Since $\beta\gamma(a) = \gamma(\beta(a)) = \gamma(a) \geq b$, it follows $\partial_{(b,c,c)}(\beta\gamma) = \beta\gamma\alpha$. Now, as in case 10.1, we obtain
$$\partial_{(b,c,c)}(\beta)\gamma  + \beta\partial_{(b,c,c)}(\gamma) =  \partial_{(b,c,c)}(\beta\gamma).$$

\emph{Case 10.5.} Let $\beta(a) \geq b$, $\gamma(a) = a$ and $\gamma(b) = c$. As in  case 10.1 we find $\partial_{(b,c,c)}(\beta) = \beta\alpha$. Since $\alpha\gamma = \overline{c}$, it follows $\partial_{(b,c,c)}(\gamma) = \overline{c}$. We obtain also $\beta\gamma = \overline{c}$ and then $\partial_{(b,c,c)}(\beta\gamma) = \overline{c}$. Now calculate
$$\partial_{(b,c,c)}(\beta)\gamma  + \beta\partial_{(b,c,c)}(\gamma) = \beta\alpha\gamma + \beta\overline{c} = \beta\overline{c} + \beta\overline{c} = \overline{c} = \partial_{(b,c,c)}(\beta\gamma).$$

\emph{Case 10.6.} Let $\beta(a) = a$, $\beta(b) = c$ and $\gamma(a) \geq b$. As in  case 10.5 we find $\alpha\beta = \overline{c}$, $\partial_{(b,c,c)}(\beta) = \overline{c}$ and $\partial_{(b,c,c)}(\gamma) = \gamma\alpha$. Then $\partial_{(b,c,c)}(\beta\gamma) = \alpha\beta\gamma + \beta\gamma\alpha = \overline{c}\gamma + \beta\gamma\alpha$ and $$\partial_{(b,c,c)}(\beta)\gamma  + \beta\partial_{(b,c,c)}(\gamma) = \overline{c}\gamma + \beta\gamma\alpha = \partial_{(b,c,c)}(\beta\gamma).$$

\emph{Case 10.7.} Let $\beta(a) = a$, $\beta(b) = b$, $\gamma(a) = a$ and $\gamma(b) = c$. As in  case 10.2 we find $\partial_{(b,c,c)}(\beta) = \beta\alpha$. As in case 10.5 we obtain $\alpha\gamma = \overline{c}$ and $\partial_{(b,c,c)}(\gamma) = \overline{c}$. Now we find $\partial_{(b,c,c)}(\beta)\gamma  + \beta\partial_{(b,c,c)}(\gamma) = \beta\alpha\gamma + \beta\overline{c} = \beta\overline{c} + \beta\overline{c} = \overline{c}$. In case 10.2 we find that $\alpha\beta(a) = b$, $\alpha\beta(b) \geq b$ and $\alpha\beta(c) \geq b$. Hence $\alpha\beta\gamma = \overline{c}$. Thus we have
$$\partial_{(b,c,c)}(\beta\gamma) = \alpha\beta\gamma + \beta\gamma\alpha = \overline{c} = \partial_{(b,c,c)}(\beta)\gamma  + \beta\partial_{(b,c,c)}(\gamma).$$

\emph{Case 10.8.} Let $\beta(a) = a$, $\beta(b) = c$, $\gamma(a) = a$ and $\gamma(b) = b$. As in the previuos case  we find $\alpha\beta = \overline{c}$ and $\partial_{(b,c,c)}(\beta) = \overline{c}$ and $\partial_{(b,c,c)}(\gamma) = \gamma\alpha$. Then, as in case 10.6 we obtain
$$\partial_{(b,c,c)}(\beta)\gamma  + \beta\partial_{(b,c,c)}(\gamma) = \overline{c}\gamma + \beta\gamma\alpha = \partial_{(b,c,c)}(\beta\gamma).$$

\emph{Case 10.9.} Let $\beta(a) = a$, $\beta(b) = c$, $\gamma(a) = a$ and $\gamma(b) = c$. As in the previuos case  we find $\alpha\beta = \overline{c}$ and $\partial_{(b,c,c)}(\beta) = \overline{c}$, $\alpha\gamma = \overline{c}$ and $\partial_{(b,c,c)}(\gamma) = \overline{c}$. Since $\beta\gamma(a) = a$ and $\beta\gamma(b) = c$, it follows $\partial_{(b,c,c)}(\beta\gamma) = \overline{c}$. Hence, we obtain
$$\partial_{(b,c,c)}(\beta)\gamma  + \beta\partial_{(b,c,c)}(\gamma) = \overline{c}\gamma + \beta\overline{c} = \overline{c} = \partial_{(b,c,c)}(\beta\gamma).$$

\emph{Case 10.10.} Let $\beta(b) \geq b$ and $\gamma \in (a,a,a)$. Easy follows that $\partial_{(b,c,c)}(\beta) = \beta\alpha = \overline{c}$. We obtain $\partial_{(b,c,c)}(\gamma) = \alpha\gamma + \gamma\alpha  = \overline{a} + \gamma\alpha  = \gamma\alpha$. Then we find $\partial_{(b,c,c)}(\beta\gamma) = \alpha\beta\gamma + \beta\gamma\alpha = \overline{a} + \beta\gamma\alpha = \beta\gamma\alpha$. So, it follows
$$\partial_{(b,c,c)}(\beta)\gamma + \beta\partial_{(b,c,c)}(\gamma) = \overline{c}\gamma + \beta\gamma\alpha = \overline{a} + \beta\gamma\alpha =  \beta\gamma\alpha = \partial_{(a,c,c)}(\beta\gamma).$$

\emph{Case 10.11.} Let $\beta \in (a,a,a)$ and $\gamma(b) \geq b$. As in the previous case $\alpha\beta = \overline{a}$ and  $\partial_{(b,c,c)}(\beta) = \beta\alpha$ and also $\partial_{(b,c,c)}(\gamma) = \gamma\alpha = \overline{c}$. Now, we find $\partial_{(b,c,c)}(\beta\gamma) = \alpha\beta\gamma + \beta\gamma\alpha = \overline{a}\gamma + \beta\gamma\alpha$.
Since $\overline{a}\gamma = \alpha\beta\gamma \leq \beta\alpha\gamma \leq \beta\gamma\alpha$, it follows $\partial_{(b,c,c)}(\beta\gamma) = \beta\gamma\alpha$. Now, we obtain
$$\partial_{(b,c,c)}(\beta)\gamma + \beta\partial_{(b,c,c)}(\gamma) = \beta\alpha\gamma + \beta\gamma\alpha =  \beta\gamma\alpha = \partial_{(b,c,c)}(\beta\gamma).$$

\emph{Case 10.12.} Let  $\beta \in (a,a,a)$ and $\gamma \in (a,a,a)$. As in the previous case $\alpha\beta = \overline{a}$ and  $\partial_{(b,c,c)}(\beta) = \beta\alpha$ and also $\alpha\gamma = \overline{a}$ and  $\partial_{(b,c,c)}(\gamma) = \gamma\alpha$. and also  $\partial_{(b,c,c)}(\beta\gamma) = \beta\gamma\alpha$. Then  we have $$\partial_{(b,c,c)}(\beta)\gamma + \beta\partial_{(b,c,c)}(\gamma) = \beta\alpha\gamma + \beta\gamma\alpha = \beta\overline{a} + \beta\gamma\alpha = \beta\gamma\alpha = \partial_{(b,c,c)}(\beta\gamma).$$

Suppose that $\beta(b) = a$, $\beta(c) \leq b$,  $\gamma(a) = a$ and $\gamma(b) = c$. Since
$$\begin{array}{c} \alpha\beta(a) = \beta(\alpha(a)) = \beta(b) = a < b = \alpha(a)  = \alpha(\beta(a)) = \beta\alpha(a),\\
\alpha\beta(b) = \beta(\alpha(b)) = \beta(c) \leq  b = \alpha(a) = \alpha(\beta(b)) = \beta\alpha(b),\\
\alpha\beta(c) = \beta(\alpha(c)) = \beta(c)  = \alpha(\beta(c)) = \beta\alpha(c),
\end{array}$$
it follows $\alpha\beta \leq \beta\alpha$ and $\partial_{(b,c,c)}(\beta) = \beta\alpha$. As in case 10.9 we obtain $\alpha\gamma = \overline{c}$ and $\partial_{(b,c,c)}(\gamma) = \overline{c}$. We find $\partial_{(b,c,c)}(\beta\gamma) = \alpha\beta\gamma + \beta\gamma\alpha$ and  $\partial_{(b,c,c)}(\beta)\gamma  + \beta\partial_{(b,c,c)}(\gamma) = \beta\alpha\gamma + \beta\overline{c} = \beta\alpha\gamma + \overline{c} = \overline{c}$. Since $\alpha(a) = b$, $\beta(b) = a$ and $\gamma(a) = a$, it follows $\alpha\beta\gamma(a) = a$. Similarly $\beta(a) = a$, $\gamma(a) = a$ and $\alpha(a) = b$ implies $\beta\gamma\alpha(a) = b$. Thus
$\partial_{(b,c,c)}(\beta\gamma)(a) = b$ and then $\partial_{(b,c,c)}(\beta\gamma) < \partial_{(b,c,c)}(\beta)\gamma  + \beta\partial_{(b,c,c)}(\gamma)$.

 From the last inequality we conclude that the endomorphisms of type $(a,a,b)$ does not belong to the subsemiring $\mathcal{D}_{(b,c,c)}$ of $\triangle^{(n)}\{a,b,c\}$, closed under the derivation $\partial_{(b,c,c)}$. If  assume that endomorphisms of type $(a,a,c)$ belongs to $\mathcal{D}_{(b,c,c)}$ we find a contradiction, because product of an endomorphism of type $(a,a,c)$ and an endomorphism of type $(a,b,b)$  is an endomorphism of type $(a,a,b)$.

 So, the semiring $\mathcal{D}_{(b,c,c)} = \mathcal{D}_{(a,c,c)}$ and does not contain endomorphisms of types $(a,a,b)$ and $(a,a,c)$. Thus we prove
\vspace{5mm}

\textbf{Proposition 10.} \textsl{The map $\partial_{(b,c,c)}$ is a derivation. Maximal subsemiring of  $\triangle^{(n)}\{a,b,c\}$ closed under this derivation is $\mathcal{D}_{(b,c,c)}$ containing all endomorphisms except the endomorphisms of types  $(a,a,b)$ and $(a,a,c)$.}
\vspace{10mm}

{\large \textbf{{3. Local derivations}}}
\vspace{4mm}

In the the all propositions from the previous section we prove for some fixed type $(x,y,z)$ of an endomorphism that $\partial_{(x,y,z)}$ is a derivation. But for different endomorphisms $\alpha_1$ and $\alpha_2$, where $\alpha_1, \alpha_2 \in (x,y,z)$ the derivations $\partial_{\alpha_1}$ and $\partial_{\alpha_2}$ are different, as we see in the following example.
\vspace{4mm}

\emph{Example 1.} In the triangle $\triangle^{(7)}\{1,3,5\}$ we consider:

--- the class $(a,a,a)$, consisting of endomorphisms $\overline{1}, 1_6 3$ and $1_6 5$,

--- the class $(a,a,b)$, consisting of endomorphisms $1_5 3_2, 1_5 3 5, 1_4 3_3$ and $1_4 3_2 5$,

--- the class $(a,a,c)$, consisting of endomorphisms $1_5 5_2$, $1_4 3 5_2$ and $1_4 5_3$.

The set of these classes is a semiring   and this semiring $S$ is closed under the derivations $\partial_{\alpha}$, where $\alpha \in S$, as we prove in cases 1,  6.1 and 9.1.

The values of derivations $\partial_{1_5 5_2}(\alpha)$, $\partial_{1_4 3 5_2}(\alpha)$ and $\partial_{1_4 5_3}(\alpha)$, where $\alpha \in S$, are shown in the following table:

$$\begin{array}{c|cccccccccc}
 & \overline{1} & 1_6 3 & 1_6 5 & 1_5 3_2 & 1_5 3 5 & 1_4 3_3 & 1_4 3_2 5 & 1_5 5_2 & 1_4 3 5_2 & 1_4 5_3\\ \hline
\partial_{1_5 5_2} & \overline{1} & \overline{1} & \overline{1} & 1_5 3_2 & 1_5 3_2 & 1_5 3_2 & 1_5 3_2 & 1_5 5_2 &1_5 5_2 & 1_4 5_3 \\
\partial_{1_4 3 5_2} & \overline{1} & \overline{1} & \overline{1} & 1_5 3_2 & 1_5 3_2 & 1_5 3_2 & 1_5 3_2 & 1_5 5_2 &1_5 5_2 & 1_4 5_3 \\
\partial_{1_4 5_3} & \overline{1} & \overline{1} & \overline{1} & 1_4 3_3 & 1_4 3_3 & 1_4 3_3 & 1_4 3_3 & 1_4 5_3 &1_4 5_3 & 1_4 5_3 \\
\end{array}.$$

Here the derivations $\partial_{1_5 5_2}$ and $\partial_{1_4 5_3}$ are different maps in the semiring $S$.

For a fixed type $(x,y,z)$ the derivations  $\partial_{\alpha}$, where $\alpha \in (x,y,z)$, are called \textbf{local derivations} of this type.

The local derivations of a given type in general case does not commute as we see in the previous example:
$$\partial_{1_4 5_3}\big(\partial_{1_5 5_2}(1_5 3_2)\big) = 1_4 3_3 > 1_5 3_3 = \partial_{1_5 5_2}\big(\partial_{1_4 5_3}(1_5 3_2)\big).$$

But for the same triangle there are commuting derivations.
\vspace{4mm}

\emph{Example 2.} Let us consider in $\triangle^{(7)}\{1,3,5\}$ the type $(a,b,b)$ which is a semiring. The endomorphisms belonging to this type are $1_3 3_4$, $1_3 3_3 5$, $1_2 3_5$ and $1_2 3_4 5$.

The values of derivations $\partial_{1_5 5_2}(\alpha)$, $\partial_{1_4 3 5_2}(\alpha)$ and $\partial_{1_4 5_3}(\alpha)$, where $\alpha$ is one of the last four endomorphisms are shown in the table:
$$\begin{array}{c|cccc}
 &  1_3 3_4 & 1_3 3_3 5 & 1_2 3_5 & 1_2 3_4 5\\ \hline
\partial_{1_5 5_2} &  1_5 3_2 & 1_5 3 5 & 1_5 3_2 & 1_5 3 5 \\
\partial_{1_4 3 5_2} & 1_4 3_3 & 1_4 3_2 5 & 1_4 3_3 & 1_4 3_2 5\\
\partial_{1_4 5_3} & 1_4 3_3 & 1_4 3_2 5 & 1_4 3_3 & 1_4 3_2 5\\
\end{array}.$$

Now, it is easy to check that derivations $\partial_{1_5 5_2}$, $\partial_{1_4 3 5_2}$ and $\partial_{1_4 5_3}$ commutes in semring of endomorphisms of type $(a,b,b)$.
\vspace{3mm}

We define addition and multiplication of types of endomorphism to be the same as operations of corresponding endomorphisms in the coordinate semiring. So, the multiplication of the types has the following table:
{\small $$
    \begin{array}{c|cccccccccc}
           \cdot \!&\! (a,a,a)\! & \!(a,a,b)\! &\! (a,a,c)\! &\! (a,b,b)\! &\! (a,b,c)\! & \!(a,c,c)\! &\! (b,b,b)\! & \! (b,b,c)\!  &\! (b,c,c)\! & \! (c,c,c)\! \\ \hline
 \!(a,a,a)\! &\! (a,a,a)\! & \!(a,a,a)\! &\! (a,a,a)\! &\! (a,a,a)\! &\! (a,a,a)\! &\! (a,a,a)\! &\! (b,b,b)\! & \!(b,b,b)\! & \! (b,b,b)\! &\! (c,c,c)\! \\
 \!(a,a,b)\! & \! (a,a,a)\! & \! (a,a,a)\! & \! (a,a,a)\! & \! (a,a,b)\! & \! (a,a,b)\! & \! (a,a,c)\! & \! (b,b,b)\! & \! (b,b,b)\! & \! (b,b,c)\! & \! (c,c,c)\! \\
 \!(a,a,c)\!& \!(a,a,a)\! & \!(a,a,b)\! & \!(a,a,c)\! & \!(a,a,b)\! & \!(a,a,c)\! & \!(a,a,c)\! & \! (b,b,b)\! & \! (b,b,c)\! & \! (b,b,c)\! & \! (c,c,c)\! \\
 \!(a,b,b)\! & \!(a,a,a)\!  & \!(a,a,a)\!  & \!(a,a,a)\!  & \!(a,b,b)\! & \!(a,b,b)\! & \!(a,c,c)\! & \! (b,b,b)\! & \! (b,b,b)\! & \! (b,c,c)\!& \! (c,c,c)\! \\
 \! (a,b,c)\!& \! (a,a,a)\! & \!(a,a,b)\! &\! (a,a,c)\! &\! (a,b,b)\! &\! (a,b,c)\! & \!(a,c,c)\! &\! (b,b,b)\! & \! (b,b,c)\!  &\! (b,c,c)\! & \! (c,c,c)\! \\
 \!(a,c,c)\! & \!(a,a,a)\! & \!(a,b,b)\! & \!(a,c,c)\! & \!(a,b,b)\! & \!(a,c,c)\! & \!(a,c,c)\! & \! (b,b,b)\! & \! (b,c,c)\! &
 \! (b,c,c)\! & \! (c,c,c)\! \\
 \! (b,b,b)\! & \! (a,a,a)\! & \! (a,a,a)\! & \! (a,a,a)\! & \! (b,b,b)\! & \! (b,b,b)\! & \! (c,c,c)\! & \! (b,b,b)\! & \! (b,b,b)\! & \! (c,c,c)\! & \! (c,c,c)\! \\
 \! (b,b,c)\! & \! (a,a,a)\! & \! (a,a,b)\! & \! (a,a,c)\! & \! (b,b,b)\! & \! (b,b,c)\! & \! (c,c,c)\! & \! (b,b,b)\! & \! (b,b,c)\! & \! (c,c,c)\! & \! (c,c,c)\! \\
 \! (b,c,c)\!& \! (a,a,a)\! & \! (a,b,b)\! & \! (a,c,c)\! & \! (b,b,b)\! & \! (b,c,c)\! & \! (c,c,c)\! &\! (b,b,b)\! & \! (b,c,c)\! & \! (c,c,c)\! & \! (c,c,c)\! \\
 \! (c,c,c)\!& \! (a,a,a)\! & \! (b,b,b)\! & \! (c,c,c)\! & \! (b,b,b)\! & \! (c,c,c)\! & \! (c,c,c)\! & \! (b,b,b)\! & \! (c,c,c)\! & \! (c,c,c)\! & \! (c,c,c)\! \\
    \end{array}.$$}

Thus this is just the table of multiplication in the semiring $\triangle^{(3)}\{0,1,2\}$.
\vspace{1mm}

Now we define a Jordan multiplication of two types $(x,y,z)$ and $(u,v,w)$ by
$$(x,y,z)\circ(u,v,w) = (x,y,z)\cdot(u,v,w) + (u,v,t)\cdot(x,y,t).$$

Then the table of the Jordan multiplication of the types of endomorphisms which are elements of $\triangle^{(n)}\{a,b,c\}$ has the following form:
{\small $$
    \begin{array}{c|cccccccccc}
           \circ \!&\! (a,a,a)\! & \!(a,a,b)\! &\! (a,a,c)\! &\! (a,b,b)\! &\! (a,b,c)\! & \!(a,c,c)\! &\! (b,b,b)\! & \! (b,b,c)\!  &\! (b,c,c)\! & \! (c,c,c)\! \\ \hline
 \!(a,a,a)\! &\! (a,a,a)\! & \!(a,a,a)\! &\! (a,a,a)\! &\! (a,a,a)\! &\! (a,a,a)\! &\! (a,a,a)\! &\! (b,b,b)\! & \!(b,b,b)\! & \! (b,b,b)\! &\! (c,c,c)\! \\
 \!(a,a,b)\! & \! (a,a,a)\! & \! (a,a,a)\! & \! (a,a,b)\! & \! (a,a,b)\! & \! (a,a,b)\! & \! (a,b,c)\! & \! (b,b,b)\! & \! (b,b,b)\! & \! (b,b,c)\! & \! (c,c,c)\! \\
 \!(a,a,c)\!& \!(a,a,a)\! & \!(a,a,b)\! & \!(a,a,c)\! & \!(a,a,b)\! & \!(a,a,c)\! & \!(a,c,c)\! & \! (b,b,b)\! & \! (b,b,c)\! & \! (b,c,c)\! & \! (c,c,c)\! \\
 \!(a,b,b)\! & \!(a,a,a)\!  & \!(a,a,b)\!  & \!(a,a,b)\!  & \!(a,b,b)\! & \!(a,b,b)\! & \!(a,c,c)\! & \! (b,b,b)\! & \! (b,b,b)\! & \! (b,c,c)\!& \! (c,c,c)\! \\
 \! (a,b,c)\!& \! (a,a,a)\! & \!(a,a,b)\! &\! (a,a,c)\! &\! (a,b,b)\! &\! (a,b,c)\! & \!(a,c,c)\! &\! (b,b,b)\! & \! (b,b,c)\!  &\! (b,c,c)\! & \! (c,c,c)\! \\
 \!(a,c,c)\! & \!(a,a,a)\! & \!(a,b,c)\! & \!(a,c,c)\! & \!(a,c,c)\! & \!(a,c,c)\! & \!(a,c,c)\! & \! (c,c,c)\! & \! (c,c,c)\! &
 \! (c,c,c)\! & \! (c,c,c)\! \\
 \! (b,b,b)\! & \! (b,b,b)\! & \! (b,b,b)\! & \! (b,b,b)\! & \! (b,b,b)\! & \! (b,b,b)\! & \! (c,c,c)\! & \! (b,b,b)\! & \! (b,b,b)\! & \! (c,c,c)\! & \! (c,c,c)\! \\
 \! (b,b,c)\! & \! (b,b,b)\! & \! (b,b,b)\! & \! (b,b,c)\! & \! (b,b,b)\! & \! (b,b,c)\! & \! (c,c,c)\! & \! (b,b,b)\! & \! (b,b,c)\! & \! (c,c,c)\! & \! (c,c,c)\! \\
 \! (b,c,c)\!& \! (b,b,b)\! & \! (b,b,c)\! & \! (b,c,c)\! & \! (b,c,c)\! & \! (b,c,c)\! & \! (c,c,c)\! &\! (c,c,c)\! & \! (c,c,c)\! & \! (c,c,c)\! & \! (c,c,c)\! \\
 \! (c,c,c)\!& \! (c,c,c)\! & \! (c,c,c)\! & \! (c,c,c)\! & \! (c,c,c)\! & \! (c,c,c)\! & \! (c,c,c)\! & \! (c,c,c)\! & \! (c,c,c)\! & \! (c,c,c)\! & \! (c,c,c)\! \\
    \end{array}.$$}

An additive subsemigrop $J$ of a semiring $S$ is said to be a Jordan ideal of $S$ if $u\circ s \in J$, for all $u \in J$ ans $s \in S$. Immediately from the last table we obtain
\vspace{3mm}

\textbf{Corollary.} \textsl{The subsemiring of $\triangle^{(n)}\{a,b,c\}$ consisting of types $(b,b,b)$, $(b,b,c)$, $(b,c,c)$ and $(c,c,c)$ is a Jordan ideal of the triangle.}
\vspace{5mm}

The last table structure gives an idea for the following result.
\vspace{3mm}

\textbf{Lemma 1.} \textsl{ Let $\alpha, \beta, \gamma \in \triangle^{(n)}\{a,b,c\}$, $\alpha$ and $\beta$ are of the same type and   $\gamma \in (a,b,c)$. Then $\gamma\alpha = \gamma\beta$.}

\emph{Proof.} $\,$ \emph{Case 1.} Let $\alpha, \beta \in (x,x,x)$, where $x = a$ or $x = b$, or $x = c$. Then $\gamma\alpha = \overline{x} = \alpha\gamma$, $\gamma\beta = \overline{x} = \beta\gamma$. So, $\gamma\alpha = \gamma\beta$.
\vspace{1mm}

\emph{Case 2.} Let $\alpha, \beta \in (a,a,b)$. Then  $\alpha = a_ib_jc_{n-i-j}$, where $a < b \leq i -1 < c \leq i + j -1$ and $\beta = a_{i_1}b_{j_1}c_{n-i_{1}-j_{1}}$, where $a < b \leq i_1 -1 < c \leq i_1 + j_1 -1$. For $\gamma = a_k b_\ell c_{n-k-\ell}$, where $a \leq k - 1< b \leq k + \ell - 1 < c$, it follows $\gamma\alpha = a_{k + \ell}b_{n-k-\ell} = \gamma\beta$.
\vspace{1mm}

\emph{Case 3.} Let $\alpha, \beta \in (a,b,b)$. Then  $\alpha = a_ib_jc_{n-i-j}$, where $a  \leq i -1 < b < c \leq i + j -1$ and $\beta = a_{i_1}b_{j_1}c_{n-i_{1}-j_{1}}$, where $a \leq i_1 -1 < b < c \leq i_1 + j_1 -1$. For $\gamma = a_k b_\ell c_{n-k-\ell}$, where $a \leq k - 1< b \leq k + \ell - 1 < c$, it follows $\gamma\alpha = a_kb_{n-k} = \gamma\beta$.
\vspace{1mm}

\emph{Case 4.} Let $\alpha, \beta \in (b,b,c)$. Then  $\alpha = a_ib_jc_{n-i-j}$, where $i \leq a < b  \leq i + j -1 < c$ and $\beta = a_{i_1}b_{j_1}c_{n-i_{1}-j_{1}}$, where $i_1 \leq a < b  \leq i_1 + j_1 -1 < c$. For $\gamma = a_k b_\ell c_{n-k-\ell}$, where $a \leq k - 1< b \leq k + \ell - 1 < c$, it follows $\gamma\alpha = b_{k + \ell}c_{n-k-\ell} = \gamma\beta$.
\vspace{1mm}

\emph{Case 5.} Let $\alpha, \beta \in (a,b,c)$. Then $\alpha$ and $\beta$ are right identities and $\gamma\alpha = \gamma\beta$.
\vspace{1mm}

\emph{Case 6.} Let $\alpha, \beta \in (a,a,c)$. Then  $\alpha = a_ib_jc_{n-i-j}$, where $ a < b\leq i-1 < i + j -1 < c$ and $\beta = a_{i_1}b_{j_1}c_{n-i_{1}-j_{1}}$, where $ a < b \leq i_{1} - 1 < i_1 + j_1 -1 < c$. For $\gamma = a_k b_\ell c_{n-k-\ell}$, where $a \leq k - 1< b \leq k + \ell - 1 < c$, it follows $\gamma\alpha = a_{k + \ell}c_{n-k-\ell} = \gamma\beta$.
\vspace{1mm}

\emph{Case 7.} Let $\alpha, \beta \in (a,c,c)$. Then  $\alpha = a_ib_jc_{n-i-j}$, where $ a \leq i-1 < i + j  \leq b < c$ and $\beta = a_{i_1}b_{j_1}c_{n-i_{1}-j_{1}}$, where $ a \leq i_{1} - 1 < i_1 + j_1 \leq b < c$. For $\gamma = a_k b_\ell c_{n-k-\ell}$, where $a \leq k - 1< b \leq k + \ell - 1 < c$, it follows $\gamma\alpha = a_{k}c_{n-k} = \gamma\beta$.
\vspace{1mm}

\emph{Case 8.} Let $\alpha, \beta \in (b,c,c)$. Then  $\alpha = a_ib_jc_{n-i-j}$, where $ i \leq a < i + j  \leq b < c$ and $\beta = a_{i_1}b_{j_1}c_{n-i_{1}-j_{1}}$, where $i_{1} \leq a < i_1 + j_1 \leq b < c$. For $\gamma = a_k b_\ell c_{n-k-\ell}$, where $a \leq k - 1< b \leq k + \ell - 1 < c$, it follows $\gamma\alpha = b_{k}c_{n-k} = \gamma\beta$.

\vspace{3mm}

We remind that the type $(a,b,c)$, that is the subsemiring  of the right identities of $\triangle^{(n)}\{a,b,c\}$ is denoted by $\mathcal{RI}\left(\triangle^{(n)}\{a,b,c\}\right)$ in [15].
\vspace{5mm}

\textbf{Theorem 1.} \textsl{The local derivations $\partial_{\alpha}$ and $\partial_{\beta}$, where $\alpha, \beta \in \triangle^{(n)}\{a,b,c\}$ are of the same type,  commutes in semiring $\mathcal{RI}\left(\triangle^{(n)}\{a,b,c\}\right)$.}

\emph{Proof.} Let us denote  $\alpha\beta = \alpha_1$ and $\beta\alpha = \beta_1$. Then $\partial_{\alpha}(\beta) = \alpha\beta + \beta\alpha  = \alpha_1 + \beta_1 = \partial_{\beta}(\alpha)$. Since $\gamma \in \mathcal{RI}\left(\triangle^{(n)}\{a,b,c\}\right)$. is a right identity, it follows $\partial_{\alpha}(\gamma) = \alpha\gamma + \gamma\alpha  = \alpha + \gamma_1$, where $\gamma_1 = \gamma\alpha$. Similarly $\partial_{\beta}(\gamma) = \beta + \gamma_2$, where $\gamma_2 = \gamma\beta$. But from Lemma 1 $\gamma_2 = \gamma_1$. So, $\partial_{\beta}(\gamma) = \beta + \gamma_1$. Now
$$\partial_{\alpha}\partial_{\beta}(\gamma) = \partial_{\beta}\big(\partial_{\alpha}(\gamma)\big) = \partial_{\beta}(\alpha + \gamma_1) = \partial_{\beta}(\alpha) + \partial_{\beta}(\gamma_1) =$$
$$= \alpha_1 + \beta_1 + \beta\gamma_1 + \gamma_1\beta = \alpha_1 + \beta_1 + \beta\gamma\alpha + \gamma\alpha\beta = \alpha_1 + \beta_1 + \beta\alpha + \gamma\alpha_1 = \alpha_1 + \beta_1 + \gamma\alpha_1.$$

 In the same way $\partial_{\beta}\partial_{\alpha}(\gamma) = \alpha_1 + \beta_1 + \gamma\beta_1$. But $\alpha_1$ and $\beta_1$ are endomorphisms of the type of $\alpha$ and $\beta$. Then from Lemma 1, it follows $\gamma\beta_1 = \gamma\alpha_1$.
 Hence, we prove $\partial_{\alpha}\partial_{\beta}(\gamma) = \partial_{\beta}\partial_{\alpha}(\gamma)$ for any $\gamma \in \mathcal{RI}\left(\triangle^{(n)}\{a,b,c\}\right)$.
\vspace{3mm}

By similar arguments one can find  local derivations which commute in another subsemiring of $\triangle^{(n)}\{a,b,c\}$.

\vspace{3mm}

\textbf{\emph{Remark.}} There is a trivial case when  local derivations commutes in $\triangle^{(n)}\{a,b,c\}$. Let $\alpha, \beta \in (c,c,c)$. For any $\gamma \in \triangle^{(n)}\{a,b,c\}$ we obtain $\partial_\alpha(\gamma) = \alpha\gamma + \gamma\alpha = \alpha\gamma + \overline{c} = \overline{c}$ and similarly $\partial_\beta(\gamma) = \overline{c}$. Since $\partial_\alpha(\partial_\beta(\gamma)) = \overline{c} = \partial_\beta(\partial_\alpha(\gamma))$, it follows $\partial_\alpha\partial_\beta = \partial_\beta\partial_\alpha$.
\newpage

{\large \textbf{{4. Jordan derivations in an arbitrary simplex}}}
\vspace{4mm}

For any endomorphism $\alpha \in \sigma^{(n)}\{a_0,\ldots,a_{k-1}\}$ we consider the map
$$\partial_\alpha : \sigma^{(n)}\{a_0,\ldots,a_{k-1}\} \rightarrow \sigma^{(n)}\{a_0,\ldots,a_{k-1}\},$$
such that $\partial_\alpha(\beta) = \alpha\beta + \beta\alpha$, where $\beta \in \triangle^{(n)}\{a,b,c\}$.

As in section 2 for any endomorphisms $\alpha, \beta, \gamma \in \sigma^{(n)}\{a_0,\ldots,a_{k-1}\}$ it follows\\
$\partial_\alpha(\beta + \gamma) = \partial_\alpha(\beta) + \partial_\alpha(\gamma)$, i.e.  the map $\partial_\alpha$ for all $\alpha \in \sigma^{(n)}\{a_0,\ldots,a_{k-1}\}$ is a linear.
\vspace{2mm}

As we know from [12] an endomorphism
 $\alpha \in \sigma^{(n)}\{a_0, a_1, \ldots, a_{k-1}\}$ is called

  \emph{endomorphism of type} $(a_{m_0}, \ldots, a_{m_{k-1}})$,  where $m_i \in \{0, \ldots, {k-1}\}$, $m_i \leq m_j$ for $i < j$, $i, j = 0, \ldots, k-1$, if $\alpha(a_i) = a_{m_i}$.

 It is clear that the relation $\alpha \sim \beta$ if and only if $\alpha$ and $\beta$ are of the same type is an equivalence relation. Any equivalence class is closed under the addition, but there are equivalence classes which are not closed under the multiplication.
\vspace{2mm}

The type of any endomorphism is itself an endomorphism of a simplex $\widehat{\mathcal{E}}_{\mathcal{C}_k} = \sigma^{(k)}\{0, 1, \ldots, k-1\}$. The simplex $\widehat{\mathcal{E}}_{\mathcal{C}_k}$ is called a \emph{coordinate simplex} of $\sigma^{(n)}\{a_0, a_1, \ldots, a_{k-1}\}$.
\vspace{3mm}

An important role in work with the types of endomorphisms is the following result
\vspace{5mm}

\textbf{Theorem 2.} [13] \emph{Let $\sigma^{(n)}\{a_0, a_1, \ldots, a_{k-1}\}$ be a simplex. Let $R$ be a subsemiring of the coordinate simplex ${\mathcal{E}}_{\mathcal{C}_k}$ of $\sigma^{(n)}\{a_0, a_1, \ldots, a_{k-1}\}$. Then the set $\widetilde{R}$ of all endomorphisms of $\sigma^{(n)}\{a_0, a_1, \ldots, a_{k-1}\}$ having a type ${\ll m_0, \ldots, m_{k-1}\gg}$, where the  endomorphism ${\wr \, m_0, \ldots, m_{k-1}\, \wr} \in R$ is a semi\-ring. Moreover, if $R$ is a (right, left) ideal of semiring $\widehat{\mathcal{E}}_{\mathcal{C}_k}$, it follows that $\widetilde{R}$ is a (right, left) ideal of simplex $\sigma^{(n)}\{a_0, a_1, \ldots, a_{k-1}\}$.}

\vspace{3mm}

As we know from section 2 the  proofs that $\partial_\alpha$\ where $\alpha \in \sigma^{(n)}\{a_0, a_1, \ldots, a_{k-1}\}$ are derivations depends of the type of this endomorphism.
\vspace{2mm}

The identity element of coordinate simplex ${\mathcal{E}}_{\mathcal{C}_k}$ forms a trivial subsemiring of ${\mathcal{E}}_{\mathcal{C}_k}$. So, from the upper theorem, it follows that endomorphisms of $\sigma^{(n)}\{a_0, a_1, \ldots, a_{k-1}\}$ of type $(0,1,\ldots,k-1)$ forms a semiring. In [13] this semiring is denoted by $\mathcal{RI}\left(\sigma^{(n)}\{a_0, a_1, \ldots, a_{k-1}\}\right)$. Elements of this semiring are right identities of the simplex and his order is $\displaystyle \prod_{i = 0}^{k-1}(a_{i+1} - a_i)$.
\vspace{3mm}

\textbf{Proposition 11.} \textsl{The map $\partial_{\alpha}$, where $\alpha \in \mathcal{RI}\left(\sigma^{(n)}\{a_0, a_1, \ldots, a_{k-1}\}\right)$ is a derivation in the whole  semiring $\sigma^{(n)}\{a_0, a_1, \ldots, a_{k-1}\}$.}

\emph{Proof.} There is not difference with the proof of Proposition 1.
\vspace{3mm}

It is easy to see that the endomorphisms of type $(a_0,a_0,\ldots,a_0)$ forms a semiring. Note that the type $(a_0,a_0,\ldots,a_0)$ is an element of the semiring $\mathcal{N}_k^{\,[0]}$ consisting of $0$--nilpotent elements of ${\mathcal{E}}_{\mathcal{C}_k}$ and order of this semiring is the $k - 1$--th Catalan number, i.e. $\displaystyle \left|\mathcal{N}_{k}^{\,[0]}\right| = \frac{1}{k}\binom{2k - 2}{k-1}$.
 \vspace{3mm}

\textbf{Proposition 12.} \textsl{The map $\partial_{\alpha}$, where $\alpha \in (a_0,a_0,\ldots,a_0)$ is a derivation in the whole semiring $\sigma^{(n)}\{a_0, a_1, \ldots, a_{k-1}\}$.}

\emph{Proof.} Exactly the same proof as for  Proposition 2.
\vspace{3mm}

 The endomorphisms of type $(a_{k-1},a_{k-1},\ldots,a_{k-1})$ also forms a semiring. This type $(a_{k-1},a_{k-1},\ldots,a_{k-1})$ is an element of the semiring $\mathcal{N}_k^{\,[k-1]}$ consisting of $k-1$--nilpotent elements of ${\mathcal{E}}_{\mathcal{C}_k}$ and order of this semiring is also the  $k - 1$--th Catalan number. \vspace{3mm}

\textbf{Proposition 13.} \textsl{The map $\partial_{\alpha}$, where $\alpha \in (a_{k-1},a_{k-1},\ldots,a_{k-1})$ is a derivation in the whole semiring $\sigma^{(n)}\{a_0, a_1, \ldots, a_{k-1}\}$.}

\emph{Proof.} Exactly the same proof as for  Proposition 3.
\vspace{1mm}

  Now we consider the endomorphisms of type $(a_{\ell},a_{\ell},\ldots,a_{\ell})$, where $0 < \ell < k-1$. They also forms a semiring. This type $(a_{\ell},a_{\ell},\ldots,a_{\ell})$ is an element of the semiring $\mathcal{N}_k^{\,[\ell]}$ consisting of $\ell$--nilpotent elements of ${\mathcal{E}}_{\mathcal{C}_k}$. We know that order of this semiring is $\left|\mathcal{N}_k^{\,[\ell]}\right| = C_\ell . C_{k-\ell-1}$, where $C_\ell$ is the $\ell$ -- th Catalan number.
   \vspace{4mm}

\textbf{Proposition 14.} \textsl{The map $\partial_{\alpha}$, where $\alpha \in (a_{\ell},a_{\ell},\ldots,a_{\ell})$ and $0 < \ell < k - 1$ is a derivation in the  semiring $\sigma^{(n)}\{a_0, a_1, \ldots, a_{k-1}\}$. Maximal subsemiring of $\sigma^{(n)}\{a_0, a_1, \ldots, a_{k-1}\}$ closed under $\partial_{\alpha}$ is the semiring $\mathcal{D}_{(a_\ell,a_\ell,\ldots,a_\ell)}$ consisting of all endomorphisms except those of types $(a_{m_0}, \ldots, a_{m_\ell}, \ldots, a_{m_{k-1}})$, where $m_\ell > \ell$ and $m_0 \leq \ell$.}

\emph{Proof.} For arbitrary endomorphisms $\alpha, \beta, \gamma \in (a_{\ell},a_{\ell},\ldots,a_{\ell})$, it follows $\beta\alpha = \overline{a_\ell}$ and then $\partial_{\alpha}(\beta) = \alpha\beta + \overline{a_\ell}$.
\vspace{2mm}

\emph{Case 1.} Let $\beta(a_\ell) \leq a_\ell$ and $\gamma(a_\ell) \leq a_\ell$. So, it follows $\alpha\beta \leq \overline{a_\ell}$, $\alpha\gamma \leq \overline{a_\ell}$ and then $\partial_{\alpha}(\beta) =  \overline{a_\ell}$ and $\partial_{\alpha}(\gamma) =  \overline{a_\ell}$. Since $\beta\gamma(a_\ell) = \gamma(\beta(a_\ell)) \leq \gamma(a_\ell) \leq a_\ell$, we obtain $\partial_{\alpha}(\beta\gamma) =  \overline{a_\ell}$. Now we find $\partial_{\alpha}(\beta)\gamma + \beta\partial_{\alpha}(\gamma) = \overline{a_\ell}\gamma + \beta\overline{a_\ell} = \overline{a_\ell}\gamma + \overline{a_\ell}$. It is clear that $\overline{a_\ell}\gamma \leq \overline{a_\ell}$. Hence $$\partial_{\alpha}(\beta)\gamma + \beta\partial_{\alpha}(\gamma) = \overline{a_\ell} = \partial_{\alpha}(\beta\gamma).$$

\emph{Case 2.} Let $\beta(a_\ell) \leq a_\ell$ and $\gamma \in (a_{j},a_{j},\ldots,a_{j})$, where $\ell < j \leq k-1$. As in the previous case $\partial_{\alpha}(\beta) =  \overline{a_\ell}$. Now, we have $\partial_{\alpha}(\gamma) = \alpha\gamma + \gamma\alpha =   \overline{a_j} + \overline{a_\ell} = \overline{a_j}$. We obtain $\partial_{\alpha}(\beta\gamma) = \alpha\beta\gamma + \beta\gamma\alpha = \overline{a_{j}} + \overline{a_\ell} = \overline{a_j}$. Thus we obtain $$\partial_{\alpha}(\beta)\gamma + \beta\partial_{\alpha}(\gamma) = \overline{a_\ell}\gamma + \beta\overline{a_j} = \overline{a_j} = \partial_{\alpha}(\beta\gamma).$$

\emph{Case 3.} Let  $\beta \in (a_j,a_j,\ldots,a_j)$, where $\ell < j \leq k-1$ and $\gamma(b) \leq b$.  As in the previous case $\partial_{\alpha}(\beta) = \overline{a_j}$ and $\partial_{\alpha}(\gamma) =  \overline{a_\ell}$. Now, it follows $\partial_{\alpha}(\beta\gamma) = \alpha\beta\gamma + \beta\gamma\alpha = \overline{a_j}\gamma + \beta\overline{a_\ell}$ and then $$\partial_{\alpha}(\beta)\gamma + \beta\partial_{\alpha}(\gamma) = \overline{a_j}\gamma + \beta\overline{a_\ell} = \partial_{\alpha}(\beta\gamma).$$

\emph{Case 4.} Let  $\beta \in (a_{k-1},a_{k-1},\ldots,a_{k-1})$ and $\gamma \in (a_{k-1},a_{k-1},\ldots,a_{k-1})$, where $\ell < j \leq k-1$. Obviously $\partial_{\alpha}(\beta) = \overline{a_j}$, $\partial_{\alpha}(\gamma) =  \overline{a_j}$, $\partial_{\alpha}(\beta\gamma) = \overline{a_j}$ and then $$\partial_{\alpha}(\beta)\gamma + \beta\partial_{\alpha}(\gamma) = \overline{a_j} = \partial_{\alpha}(\beta\gamma).$$

Suppose that $\beta(a_\ell) = a_0$, $\gamma(a_\ell) = a_j$, where $\ell < j \leq k-1$ and $\gamma(a_0) = a_i$, where $i < j$. Now $\alpha\beta = \overline{a_0}$ and $\beta\alpha = \overline{a_\ell}$, so, $\partial_{\alpha}(\beta) = \overline{a_\ell}$. We obtain $\gamma\alpha = \overline{a_\ell}$ and $\alpha\gamma = \overline{a_j}$ and then $\partial_{\alpha}(\gamma) = \overline{a_j}$. Now we find $\partial_{\alpha}(\beta\gamma) = \alpha\beta\gamma + \beta\gamma\alpha = \overline{a_0}\gamma + \beta\overline{a_\ell} = \overline{a_0}\gamma + \overline{a_\ell}$. But $\overline{a_0}\gamma = \overline{a_i}$. If $s = \max\{i,\ell\}$, then $\partial_{\alpha}(\beta\gamma) = \overline{a_s}$. Now, it follows
$$\partial_{\alpha}(\beta)\gamma + \beta\partial_{\alpha}(\gamma) = \overline{a_\ell}\gamma + \beta\overline{a_j} = \overline{a_j}  + \overline{a_j} = \overline{a_j} > \overline{a_s} = \partial_{\alpha}(\beta\gamma).$$

From the last inequality foolows that  the endomorphisms of all types $(a_{m_0}, \ldots, a_{m_\ell}, \ldots, a_{m_{k-1}})$, where $m_\ell > \ell$ and $m_0 \leq \ell$ does not belong to the subsemiring $\mathcal{D}_{(a_\ell,a_\ell,\ldots,a_\ell)}$, closed under derivation   $\partial_{\alpha}$, where $\alpha \in (a_{k-1},a_{k-1},\ldots,a_{k-1})$.
\vspace{3mm}

For a fixed type $(a_{m_0},  \ldots, a_{m_{k-1}})$ the derivations  $\partial_{\alpha}$, where\\ $\alpha \in (a_{m_0}, \ldots, a_{m_\ell}, \ldots, a_{m_{k-1}})$, are called \textbf{local derivations} of this type.
\vspace{3mm}

Now, let endomorphism $\alpha$ be of type $(a_{m_0}, \ldots, a_{m_{k-1}})$ and the endo\-morphism $(m_0, \ldots, m_{k-1})$ from the coordinate simplex ${\mathcal{E}}_{\mathcal{C}_k}$ be an idem\-potent, different from all the constant endo\-mor\-phisms $\overline{j}$, where $j = 0, \ldots, k -1$, and from the identity $\mathbf{i}$. Then we say that $\alpha$ is an endomorphism of an \emph{idempotent form}.
\vspace{2mm}

\textbf{Lemma 2.} \textsl{ Let:}

\textsl{ a)} $\alpha, \beta \in \mathcal{RI}\left(\sigma^{(n)}\{a_0, a_1, \ldots, a_{k-1}\}\right)$;

\textsl{ b)} $\alpha, \beta \in (a_\ell, a_\ell, \ldots, a_\ell)$, \textsl{where} $\ell = 0, \ldots, k-1$;

  \textsl{ c)}  $\alpha$ and $\beta$ \textsl{are of the same type and also of an idempotent form.}

  \textsl{Let} $\gamma \in \mathcal{RI}\left(\sigma^{(n)}\{a_0, a_1, \ldots, a_{k-1}\}\right)$. \textsl{Then}
   $\gamma\alpha = \gamma\beta$.

\vspace{2mm}

\emph{Proof.}  Let $\alpha = (a_{m_0})_{i_0}(a_{m_1})_{i_1} \cdots (a_{m_{k-1}})_{i_{k-1}}$ and $\beta = (a_{m_0})_{j_0}(a_{m_1})_{j_1} \cdots (a_{m_{k-1}})_{j_{k-1}}$.

 a) Since $\alpha$  and $\beta$ are right identities then $\gamma\alpha = \gamma = \gamma\beta$.

 b) Now $\gamma\alpha = \overline{a_\ell} = \alpha\gamma$, $\gamma\beta = \overline{a_\ell} = \beta\gamma$. So, $\gamma\alpha = \gamma\beta$.

 c) Let $\gamma = (a_0)_{\ell_0}(a_1)_{\ell_1} \cdots (a_{k-1})_{\ell_{k-1}}$ and $\ell_{m_s}$ is the greatest integer which is less or equal to the number $m_s$, where $s = 0, 1, \ldots, k-1$. Then from (3) of 6.1, it follows
$$\gamma\alpha = (a_{m_0})_{\sum_{p=0}^{m_0} \ell_p}(a_{m_1})_{\sum_{p=m_0 + 1}^{m_1} \ell_p} \cdots (a_{m_{k-1}})_{\sum_{p=m_{k-2} + 1}^{m_{k-1}} \ell_p} = \gamma\beta.$$
\vspace{1mm}

\textbf{Theorem 3.} \textsl{The local derivations $\partial_{\alpha}$ and $\partial_{\beta}$, where $\alpha$ and $\beta$ are endomorphisms of the same type, such that:}

\textsl{ a)} $\alpha, \beta \in \mathcal{RI}\left(\sigma^{(n)}\{a_0, a_1, \ldots, a_{k-1}\}\right)$;

\textsl{ b)} $\alpha, \beta \in (a_\ell, a_\ell, \ldots, a_\ell)$, \textsl{where} $\ell = 0, \ldots, k-1$;

  \textsl{ c)}  $\alpha$ and $\beta$ \textsl{are of the same type and also of an idempotent form,}

\textsl{
commutes in semiring $\mathcal{RI}\left(\sigma^{(n)}\{a_0, a_1, \ldots, a_{k-1}\}\right)$.}

\emph{Proof.} Exactly the same proof as for  Theorem 1.

\vspace{3mm}

\textbf{\emph{Remark.}} The number of all types considered in the last theorem is the even Fibonacci number $f_{2k}$  -- [3], Theorem 14.3.6.

\vspace{8mm}

{\large \textbf{{5. Jordan derivation which is not a derivation}}}
\vspace{4mm}

Here we use some notations and facts from [20].
\vspace{2mm}

  Let us consider the upper triangular Toeplitz $n\times n$ matrices with entries from additively idempotent semiring $S_0$.
In [20] the additively idempotent semiring of such matrices we  denote by $\mathbb{UT}_n(S_0)$.
\vspace{2mm}

Remain that any matrix $A \in \mathbb{UT}_n(S_0)$ has a form
$$A = a_0 E + a_1 D + \cdots + a_{n-1} D^{n-1},$$
where $a_i \in S_0$, $i = 0, 1, \ldots, n-1$, $E$ is an identity matrix and
$$D = \left(
    \begin{array}{cccccc}
      0 & 1 & 0 & \cdots & 0 & 0 \\ \vspace{1mm}
      0 & 0 & 1 &\cdots & 0  & 0\\ \vspace{1mm}
      \vdots & \vdots& \vdots & \vdots & \vdots & \vdots \\
      0 & 0 & 0& \cdots & 0 & 1 \\
      0 & 0 & 0 & \cdots & 0 & 0 \\
    \end{array}
  \right).$$

  Let $X \in \mathbb{UT}_n(S_0)$ be a fixed matrix. For any $A \in \mathbb{UT}_n(S_0)$ we have
  $$AX = a_0x_0E + (a_0x_1 + a_1x_0)D + \cdots + \sum_{i=0}^k a_ix_{k-i} D^k + \cdots + (a_0x_{n-1} + a_{n-1}x_0)D^{n-1}$$
  and
  $$XA = x_0a_0E + (x_0a_1 + x_1a_0)D + \cdots + \sum_{i=0}^k x_ia_{k-i} D^k + \cdots + (x_0a_{n-1} + x_{n-1}a_0)D^{n-1}.$$

  Now, the Jordan multiplication $a\circ b = ab + ba$ in semiring $S_0$ induces a Jordan multiplication in $\mathbb{UT}_n(S_0)$.
  Thus we have
  $$X\circ A = AX + XA = (x_0\circ a_0)E + (x_0\circ a_1 + x_1\circ a_0)D + \cdots +$$ $$+ \sum_{i=0}^k (x_i\circ a_{k-i}) D^k + \cdots + (x_0\circ a_{n-1} + x_{n-1}\circ a_0)D^{n-1} \in \mathbb{UT}_n(S_0).$$

  Let us define $\delta_X(A) = X\circ A$. Then for any $B \in \mathbb{UT}_n(S_0)$, it follows $\delta_X(B) = X\circ B$ and $\delta_X(A\circ B) = X\circ A\circ B$. Since $\delta_X(A)\circ B = X\circ A\circ B$ and $A\circ \delta_X(B) = A\circ X\circ B = (A\circ X)\circ B = (X\circ A)\circ B = X\circ A\circ B$, it follows
  $$\delta_X(A\circ B) = \delta_X(A)\circ B + A\circ \delta_X(B).$$

  Thus we prove that $\delta_X$ is a Jordan derivation.
  \vspace{2mm}

  On the other hand we consider $\delta_X(A\, B) = X\circ AB = XAB + ABX$.

  We calculate $\delta_X(A)\, B = (XA + AX)B = XAB + AXB$ and $A\,\delta_X(B) = A(XB + BX) = AXB + ABX$. In the general case
  $$\delta_X(A\, B) = XAB + ABX \neq XAB + AXB + ABX = \delta_X(A)\, B  + A\,\delta_X(B),$$
  as we understand by the following resoning:
  \vspace{2mm}

  The element $a_{11}$ in the matrix in the left side of the inequality above is equal to $$x_0a_0b_0 + a_0b_0x_0,$$ but this element in the matrix in the right side is equal to $$x_0a_0b_0 + a_0x_0b_0 + a_0b_0x_0.$$

  Hence $\delta_X$ is a Jordan derivation, but is not a derivation.
\newpage

\vspace{8mm}

{\large  \textbf{References}}

\vspace{3mm}

[1] Bre\v{s}ar M., Vukman J.,  Jordan derivations on prime rings, Bull. Austral. Math. Soc. 37, 1988, 321-322.
\vspace{1mm}

[2] Bre\v{s}ar M.,  Jordan derivations on semiprime rings, Proc. Amer. Math. Soc. 104 (4), 1988, 1003-1006.
\vspace{1mm}

 [3] Ganyushkin O., Mazorchuk V., Classical Finite Transformation Semigroups: An Introduction,
Springer-Verlag London Limited, 2009.
\vspace{1mm}

 [4] Golan J., Semirings and Affine Equations over Them: Theory and Applications, Springer, 2003.
\vspace{1mm}

[5]   Herstein I. N., Jordan derivations of prime rings, Proc. Amer. Math. Soc. 8, 1957, 1104-1110.
\vspace{1mm}

[6] Herstein I. N.,  Lie and Jordan structures in simple, associative rings, Bull. Amer. Math. Soc. 67 (6), 1961
517-531.
\vspace{1mm}

  [7] Kolchin E. R., Differential Algebra and Algebraic Groups, Academic Press, New
  York,  London, 1973.
\vspace{1mm}

[8] Kozlov D., Combinatorial Algebraic Topology, Springer--Verlag Berlin, 2008.
\vspace{1mm}

[9] Ritt J. F., Differential Algebra, Amer. Math. Soc. Colloq. Publ. 33, New York, 1950.
\vspace{1mm}

[10] Stanley R., Enumerative combinatorics, Vol. 2, Cambr, Univ. Press, 1999.
\vspace{1mm}

[11] Trendafilov I., {Derivations in Some Finite Endomorphism Semirings}, Discuss.   Math.
Gen. Algebra and Appl., Vol. 32, 2012, 77--100.
\vspace{1mm}

 [12]  Trendafilov I.,  A partition of geometrical structures of the endomorphism semiring, Comptes rendus de l'Acad\'{e}mie bulgare des Sciences, Tome 66, No 12, 2013, 1661 -- 1668.
\vspace{1mm}

[13] Trendafilov I., {Simplices in the Endomorphism Semiring of a Finite Chain},  Hindawi Publishing Corporation,
Algebra, Volume 2014, Art. ID 263605, 2014.
\vspace{1mm}

[14] Trendafilov I., Vladeva D., {The endomorphism semiring of a finite chain}, Proc.  Techn. Univ.-Sofia, 61, 1, 2011, 9 -- 18.
\vspace{1mm}

[15] Trendafilov I., Vladeva D., Combinatorial Results for Geometric Structures in   Endomorphism Semirings,   AIP Conf. Proc.  1570, 2013,  461 -- 468.
\vspace{1mm}
\vspace{1mm}

[16] Trendafilov I., Vladeva D., Derivations in a triangle -- I part. The projection on the least string of a triangle is a derivation, Proc. Techn. Univ.-Sofia, 65, 1, 2015, 243 -- 252.
\vspace{1mm}

[17] Trendafilov I.,  Vladeva D., Derivations in a triangle -- II part. The projection  on any of strings of a triangle is a derivation, Proc. Techn. Univ.-Sofia, 65, 1, 2015,  253 -- 262
\vspace{1mm}

[18]  Vladeva D., Derivations in a endomorphism semiring,
Serdica Math. J.-- Bulgarian Academy of Sciences, Vol.42, No. 3-4, pp. 251 -- 260.

[19] Vladeva D., Projections of $k$--simplex onto the subsimplices of arbitrary type are derivations, arXiv:1706.00033v1 [math.RA] 31 May 2017.
\vspace{1mm}

[20] Vladeva D., Derivatives of triangular, Toeplitz, circulant matrices
and of matrices of other forms over semirings, arXiv:1707.04716v1 [math.RA] 15 Jul 2017.
\vspace{7mm}

\begin{flushleft} { Dimitrinka I. Vladeva}\\
Department ``Mathema\-tics and Physics''\\ University of Forestry, Sofia, Bulgaria,\\ \emph{e-mail:}
d$\_$vladeva@abv.bg
\end{flushleft}

\end{document}